\documentclass[12pt]{amsart} \textwidth=14.5cm \oddsidemargin=1cm
\usepackage{amsmath}
\usepackage{amsxtra}
\usepackage{amscd}
\usepackage{amsthm}
\usepackage{amsfonts}
\usepackage{amssymb}
\usepackage{eucal}

\input prepictex
\input pictex
\input postpictex

\newtheorem{thm}{Theorem}[section]
\newtheorem{lem}{Lemma}[section]
\newtheorem{cor}{Corollary}[section]
\newtheorem{prop}{Proposition}[section]

\theoremstyle{definition}

\theoremstyle{remark}
\newtheorem{rem}{Remark}[section]

\numberwithin{equation}{section}

\begin{document}

\newcommand{\thmref}[1]{Theorem~\ref{#1}}
\newcommand{\secref}[1]{Section~\ref{#1}}
\newcommand{\lemref}[1]{Lemma~\ref{#1}}
\newcommand{\propref}[1]{Proposition~\ref{#1}}
\newcommand{\corref}[1]{Corollary~\ref{#1}}
\newcommand{\remref}[1]{Remark~\ref{#1}}
\newcommand{\eqnref}[1]{(\ref{#1})}
\newcommand{\exref}[1]{Example~\ref{#1}}

\newcommand{\nc}{\newcommand}
\nc{\on}{\operatorname} \nc{\Z}{{\mathbb Z}} \nc{\C}{{\mathbb C}}
\nc{\R}{{\mathbb R}} \nc{\oo}{{\mf O}} \nc{\N}{{\mathbb N}}
\nc{\bib}{\bibitem} \nc{\pa}{\partial} \nc{\F}{{\mf F}}
\nc{\rarr}{\rightarrow} \nc{\larr}{\longrightarrow}
\nc{\al}{\alpha} \nc{\ri}{\rangle} \nc{\lef}{\langle} \nc{\W}{{\mc
W}} \nc{\gam}{\ol{\gamma}} \nc{\Q}{\ol{Q}} \nc{\q}{\widetilde{Q}}
\nc{\la}{\lambda} \nc{\ep}{\epsilon} \nc{\g}{\mf g} \nc{\h}{\mf h}
\nc{\n}{\mf n} \nc{\bb}{\mf b} \nc{\A}{{\mf a}} \nc{\G}{{\mf g}}
\nc{\D}{\mc D} \nc{\Li}{{\mc L}} \nc{\La}{\Lambda}
\nc{\is}{{\mathbf i}} \nc{\V}{\mf V} \nc{\bi}{\bibitem}
\nc{\NS}{\mf N} \nc{\dt}{\mathord{\hbox{${\frac{d}{d t}}$}}}
\nc{\E}{\mc E} \nc{\ba}{\tilde{\pa}} \nc{\half}{\frac{1}{2}}
\def\smapdown#1{\big\downarrow\rlap{$\vcenter{\hbox{$\scriptstyle#1$}}$}}
\nc{\mc}{\mathcal} \nc{\mf}{\mathfrak} \nc{\ol}{\fracline}
\nc{\el}{\ell} \nc{\etabf}{{\bf \eta}} \nc{\x}{{\bf x}}
\nc{\xibf}{{\bf \xi}} \nc{\y}{{\bf y}} \nc{\WW}{\mc W}
\nc{\SW}{\mc S \mc W} \nc{\sd}{\mc S \mc D} \nc{\hsd}{\widehat{\mc
S\mc D}} \nc{\parth}{\partial_{\theta}} \nc{\cwo}{\C[w]^{(1)}}
\nc{\cwe}{\C[w]^{(0)}} \nc{\hf}{\frac{1}{2}}
\nc{\hsdzero}{{}^0\widehat{\sd}} \nc{\hsdpp}{{}^{++}\widehat{\sd}}
\nc{\hsdpm}{{}^{+-}\widehat{\sd}}
\nc{\hsdmp}{{}^{-+}\widehat{\sd}}
\nc{\hsdmm}{{}^{--}\widehat{\sd}} \nc{\gltwo}{{\rm
gl}_{\infty|\infty}} \nc{\btwo}{{B}_{\infty|\infty}}
\nc{\htwo}{{\h}_{\infty|\infty}} \nc{\hglone}{\widehat{\rm
gl}_{\infty}} \nc{\hgltwo}{\widehat{\rm gl}_{\infty|\infty}}
\nc{\hbtwo}{\hat{B}_{\infty|\infty}}
\nc{\hhtwo}{\hat{\h}_{\infty|\infty}} \nc{\glone}{{\rm gl}_\infty}
\nc{\gl}{{\rm gl}} \nc{\ospd}{\mc B} \nc{\hospd}{\widehat{\mc B}}
\nc{\pd}{\mc P} \nc{\hpd}{\widehat{\pd}} \nc{\co}{\mc O}
\nc{\Oe}{\co^{(0)}} \nc{\Oo}{\co^{(1)}} \nc{\sdzero}{{}^0{\sd}}
\nc{\hz}{\hf+\Z} \nc{\vac}{|0 \rangle} \nc{\K}{\mf k}
\nc{\bhf}{\bf\hf}

\advance\headheight by 2pt

\title[Howe Duality and Combinatorial Character Formula]
{Howe Duality and Combinatorial Character Formula for
Orthosymplectic Lie superalgebras}

\author[Shun-Jen Cheng]{}
\thanks{$^*$Partially supported by NSC-grant 91-2115-M-002-007 of the R.O.C}

\author[R.~B.~Zhang]{}

\maketitle

\centerline{Shun-Jen Cheng${}^*$}
\medskip\centerline{\Small{Department of Mathematics, National
Taiwan University}} \centerline{\Small{Taipei, Taiwan 106
}}\centerline{\Small{\em E-mail: chengsj@math.ntu.edu.tw}}

\bigskip
\bigskip

\centerline{R.~B.~Zhang}\medskip \centerline{\Small{School of
Mathematics and Statistics, University of Sydney}}
\centerline{\Small{New South Wales 2006, Australia}}
\centerline{\Small{\em E-mail: rzhang@maths.usyd.edu.au}}

\begin{abstract}
We study the Howe dualities involving the reductive dual pairs
$(O(d),spo(2m|2n))$ and $(Sp(d),osp(2m|2n))$ on the
(super)symmetric tensor of $\C^d\otimes\C^{m|n}$. We obtain
complete decompositions of this space with respect to their
respective joint actions.  We also use these dualities to derive a
character formula for these irreducible representations of
$spo(2m|2n)$ and $osp(2m|2n)$ that appear in these decompositions.
\end{abstract}

\tableofcontents

\section{Introduction}

Howe duality \cite{H1,H2} relates the representation theories of a
pair of Lie groups/ algebras. It enables the study of
representations of one Lie group/algebra via the representations
of its dual partner, and hence it has become a fundamental tool
where representation theory of classical Lie groups/algebras is
indispensable. As simple and fundamental a concept it is therefore
of no surprise that the Howe duality also applies to
generalizations of finite-dimensional Lie groups/algebras.  We
point out here the Howe dualities of finite-dimensional Lie
superalgebras in \cite{CW1,CW2,LS,N2,S1,S2} of
infinite-dimensional Lie algebras in \cite{W1,W2,KWY,H,FKRW} and
of infinite-dimensional Lie superalgebras in \cite{CW3}. In the
above-mentioned articles, the main themes revolve around the
construction of Howe dualities. In the present article we are also
concerned about applications of the Howe dualities that we obtain.

Consider a Lie superalgebra whose representation theory we wish to
study. Suppose that on some natural space one has a Howe duality
involving this Lie superalgebra with a classical Lie group or Lie
algebra as its dual partner.  As the representation theory of its
classical counterpart is well-understood, one expects that this
should enable one to study the representations of the Lie
superalgebra in question with the help of the representation
theory of its classical dual partner. Of particular interest is a
derivation of a character formula for this Lie superalgebra.  It
appears plausible that knowing the character of the total space
and the characters of each of the irreducible representations of
the classical group/algebra, one should in principle be able to
obtain a character formula for the Lie superalgebra in question.
As is well-known, character formulas for Lie superalgebras in
general are rather difficult to obtain, and hence such a method
could facilitate the computation of characters for certain
representations of Lie superalgebras. One of the main purposes of
this paper is to demonstrate for the orthosymplectic Lie
superalgebra that such an approach to character formulas is indeed
viable. The general idea is the following.

Let $\G_m$ be a classical Lie algebra of rank $m$ and let $X$ be a
fixed finite-dimensional classical Lie algebra. Suppose on some
space $\F_m$ the pair $(\G_m,X)$ forms a dual pair in the sense of
Howe. Suppose that this is the case for every $m$. That is, we
have for each $m$ a (multiplicity-free) decomposition with respect
to $\G_m\times X$ of the form
\begin{equation*}
\F_m=\sum_{\la}V_{\G_m}^{\la}\otimes V_X^{\la'},
\end{equation*}
where $V_{\G_m}^\la$ and $V_X^{\la'}$ denote irreducible
representations of $\G_m$ and $X$, respectively.  Here $\la$ is
summed over a subset of irreducible representations of $\G_m$.
Since here the correspondence between irreducible representations
of $\G_m$ and $X$, given by $\la\rightarrow\la'$, is one-to-one,
we will write $V_X^\la$ for $V_X^{\la'}$.

Now suppose that $\G_{m|n}$ is the Lie superalgebraic analogue of
$\G_m$ and we have an action of the dual pair $\G_{m|n}\times X$
on $\F_{m|n}$, which is the tensor product of $\F_m$ with a
Grassmann superalgebra depending on $n$. Thus we have similarly
\begin{equation*}
\F_{m|n}=\sum_{\la}V_{\G_{m|n}}^\la\otimes V_X^\la,
\end{equation*}
where $V_{\G_{m|n}}^\la$ denotes an irreducible representation of
$\G_{m|n}$.

Our claim is that if one knows the characters of $V_{\G_m}^\la$
for every $m$, then one, in principle, also knows the characters
for $V^\la_{\G_{m|n}}$.

Let us now discuss the content of the present article in more
detail. Let $X=O$ or $X=Sp$ so that $X(d)$ denotes either the
orthogonal or the symplectic group acting on $\C^d$. We have an
induced action on $\C^d\otimes\C^m$, thus giving rise to an action
on the symmetric tensor $S(\C^d\otimes\C^m)$. Now by classical
invariant theory (c.f.~\cite{H1,GW}) the invariants of $X(d)$ in
the endomorphism ring of $S(\C^d\otimes\C^m)$ is generated by
quadratic invariants, which may be identified with the Lie algebra
$sp(2m)$ in the case $X=O$ and $so(2m)$ in the case $X=Sp$. This
implies that $(O(d),sp(2m))$ and $(Sp(d),so(2m))$ are Howe dual
pairs on $S(\C^d\otimes\C^m)$ .

Now let $\C^{m|n}$ be the complex superspace of dimension $(m|n)$.
The Lie group $X(d)$ acts in a similar fashion on the (super)
symmetric tensor $S(\C^d\otimes\C^{m|n})$.  Analogously one
derives the $(O(d),spo(2m|2n))$- and the $(O(d),osp(2m|2n))$-Howe
duality on $S(\C^d\otimes\C^{m|n})$. Although these dualities
appear already in Howe's classical paper \cite{H1}, the complete
decompositions of $S(\C^d\otimes\C^{m|n})$ with respect to these
joint actions are unknown to the best of our  knowledge. In
\cite{N2} a partial decomposition is obtained for $X=O$, with a
complete answer given in the case of $m=n=1$ only.

Our first main task is to give the complete decompositions of
$S(\C^d\otimes\C^{m|n})$ with respect to these Howe dual pairs.
This is achieved in the following way. By \cite{H1} the
decomposition of $S(\C^d\otimes\C^{m|n})$ with respect to $X(d)$
and its dual partner is reduced to the decomposition of the
subspace of harmonic polynomials $H$ with respect to the dual pair
$(X(d),gl(m|n))$. Our task is then reduced to the construction of
all $(X(d),gl(m|n))$-highest weight vectors in $H$. Our analysis
of the $(X(d),gl(m|n))$-highest weight vectors in $H$ relies
heavily on the $(gl(d),gl(m|n))$-Howe duality in
$S(\C^d\otimes\C^{m|n})$ in \cite{CW1} (see also \cite{S1,S2}) and
the description of their joint highest weight vectors given in
\cite{CW1}. Another important ingredient is the construction of an
explicit basis for each irreducible $gl(d)\times
gl(m|n)$-component that appears in $S(\C^d\otimes\C^{m|n})$.

The idea to obtain a character formula for the irreducible
representations of $spo(2m|2n)$ or $osp(2m|2n)$ is roughly as
follows. In order to simplify notation we take $X=Sp$ in what
follows, but note that the same applies to $X=O$ with minor
modification. We first consider the classical duality, i.e.~the
case when $n=0$. Thus we have an identity of characters of the
form
\begin{equation*}
{\rm ch}S(\C^d\otimes\C^m)=\sum_{\la}{\rm
ch}V^\la_{Sp(d)}\otimes{\rm ch}V^\la_{so(2m)}.
\end{equation*}
Since now characters are polynomial functions on a Cartan
subalgebra, we can write $\chi_{Sp(d)}^\la(\x)={\rm
ch}V^\la_{Sp(d)}$ and $\chi_{so(2m)}^\la(\y)={\rm
ch}V^\la_{so(2m)}$, where $\x$ and $\y$ denote the linear
functions on the respective Cartan subalgebras. The left-hand side
is the character of the algebra of polynomials in $dm$ variables,
which is a symmetric function in $\x$ and $\y$. Now taking the
limit as $m\to\infty$ in an appropriate way one obtains a
combinatorial identity involving infinitely many variables
$\y=y_1,y_2,\cdots,y_m,\cdots$. Since the right-hand side is
symmetric in $\y$, we may apply to this identity the involution
$\omega$ of symmetric functions that sends the complete symmetric
functions to the elementary symmetric functions (see \cite{M}).
The $\omega$ turns the left-hand side into the character of the
tensor product of a polynomial algebra with a Grassmann algebra.
Therefore, due to ``linear independence'' of the
$\chi^\la_{Sp(d)}$, it follows that (modulo some minor
manipulation of the variables) the expression
$\omega(\lim_{m\to\infty}\chi^\la_{so(2m)}(\y))$ is essentially
the character of the irreducible representation of $osp(2m|2n)$
paired with $V_{Sp(d)}^\la$.  At this point we wish to point out
our results imply that the characters of the representations of
the Lie superalgebra $osp(2m|2n)$ (respectively $spo(2m|2n)$), for
any $m,n\in\Z_+$, that appear under the Howe duality are
completely determined by the characters of the representations of
the Kac-Moody algebra corresponding to the infinite affine matrix
$D_\infty$ (respectively $C_\infty$) (see \cite{K2}) that appear
under a similar Howe duality.

The next problem is to describe the expression
$\chi^\la_{so(2m)}(\y)$. For this we use the beautiful formula of
Enright \cite{E,DES} for unitarizable irreducible representations
associated to a classical Hermitian symmetric pair.  The reason
for this is that in our case we may express such a character in
terms of Schur functions which are carried by $\omega$ to the
so-called hook Schur functions of Berele and Regev \cite{BR}. This
allows us to obtain a satisfactory description of the characters.

We now come to the organization of the paper. In
\secref{parameterization} we recall some basic facts on the
orthogonal and symplectic groups and the orthosymplectic
superalgebra, where we also take the opportunity to set the
notation to be used throughout the paper. In \secref{glglduality}
we recall the $(gl(d),gl(m|n))$-duality on
$S(\C^d\otimes\C^{m|n})$ and construct an explicit basis for each
irreducible component that appears in the decomposition of
$S(\C^d\otimes\C^{m|n})$. In \secref{duality} and
\secref{jointhwv} we study the $(O(d),spo(2m|2n))$-duality and the
$(Sp(d),osp(2m|2n))$-duality and obtain the complete
decompositions of $S(\C^d\otimes\C^{m|n})$ with respect to their
respective joint actions. In \secref{character} we derive a
character formula for these representations of $spo(2m|2n)$ and
$osp(2m|2n)$ that appear in the decomposition of
$S(\C^d\otimes\C^{m|n})$.  Here we should mention that in the case
of $O(d)\times spo(2m|2n)$ with $d$ even, we are only able to
derive the formula for a sum of two irreducible representations in
general. We also remark that in \cite{N1} a character formula for
the oscillator representations is given. This corresponds to our
case $O(1)$. In order to obtain a better description of the
character formulas we are required to study Enright's formula in
more detail.  This is done in the \secref{group}. In
\secref{consequences} we study the character formulas in more
detail. In \secref{tensor}, as another application of our Howe
dualities, we give formulas for decomposing tensor products of
these irreducible $spo(2m|2n)$- and $osp(2m|2n)$ modules that
appear in the decomposition of $S(\C^d\otimes\C^{m|n})$.

Finally all vector spaces, algebras etc.~are over the complex
field $\C$ unless otherwise specified.  By a {\em partition} we
mean a non-increasing finite sequence of non-negative integers. A
{\em composition} is a finite sequence of either all non-negative
integers or all positive half-integers. Furthermore, by a {\em
generalized partition} we will always mean a finite non-increasing
sequence of either all integers or all half integers.  By a {\em
generalized composition} we will mean a finite sequence of either
all integers or all half integers.

\section{Parameterization of irreducible
representations}\label{parameterization}

In this section we give parameterizations of irreducible
representations of the Lie groups and Lie superalgebras that we
will be dealing with in this paper. For a more complete treatment
of the material on Lie groups the reader is referred to \cite{BT}.

\subsection{Irreducible representations of the general linear Lie
superalgebra}\label{irrepglmn} Let $\C^{m|n}$ denote the complex
$(m|n)$-dimensional superspace. The space of complex linear
transformations on $\C^{m|n}$ has a natural structure as a Lie
superalgebra, which we will denote by $gl(m|n)$.  Choose a
homogeneous basis for $\C^{m|n}$ so that we may regard $gl(m|n)$
as $(m+n)\times(m+n)$ matrices.  Denote by $E_{ij}$ the elementary
matrix with $1$ in the $i$-th row and $j$-th column and $0$
elsewhere.  Then $\h=\sum_{i}\C E_{ii}$ is a Cartan subalgebra,
while $B=\sum_{i\le j}\C E_{ij}$ is a Borel subalgebra containing
$\h$.  Recall that finite-dimensional irreducible
$gl(m|n)$-modules are parameterized by $\la\in\h^*$ with
$\la_i-\la_{i+1}\in\Z_+$, for $i=1,\cdots,m-1, m+1,\cdots,m+n-1$,
where $\la_i=\la(E_{ii})$. We will denote the corresponding
finite-dimensional irreducible module by $V^\la_{m|n}$. Suppose
that $\la$ is a partition (or a Young diagram) with $\la_{m+1}\le
n$. Then drawing the corresponding diagram $\la$ may be visualized
as lying in the $(m|n)$-hook, i.e.~from $n+1$-st column on the
columns of $\la$ all have lengths less than $m+1$. We may
interpret $\la$ as a highest weight of $gl(m|n)$ by associating to
the diagram $\la$ the labels
$\la=(\la_1,\cdots,\la_m;<\la'_1-m>,\cdots,<\la'_n-m>)$, where
$\la'_i$ is the length of the $i$-th column of the diagram $\la$,
and $<r>$ stands for $r$, if $r\in\N$, and $0$ otherwise. If clear
from the context that $\la$ is a Young diagram with $\la_{m+1}\le
n$, we will mean by $V^\la_{m|n}$ the irreducible $gl(m|n)$-module
of highest weight $\la$.

\subsection{Irreducible representations of the orthogonal
group} Let us denote by $\{e^1,\cdots,e^d\}$ the standard basis
for $\C^d$. Consider the symmetric non-degenerate bilinear form
determined by the $d\times d$ matrix
\begin{equation*}
J_d=\begin{pmatrix}
0&0&\cdots&0&1\\
0&0&\cdots&1&0\\
\vdots&\vdots&\vdots&\vdots&\vdots\\
0&1&\cdots&0&0\\
1&0&\cdots&0&0\\
\end{pmatrix}.
\end{equation*}
The complex orthogonal group $O(d)$ is the subgroup of the complex
general linear group $GL(d)$ preserving this form.  The Lie
algebra of $O(d)$ is $so(d)$, which consists of those $A\in gl(d)$
with $J_dA^tJ_d+A=0$, that is, $A$ is skew-symmetric with respect
to the diagonal running from the top right to the bottom left
corner.

Consider the case when $d=2k$ is even. We take as a Borel
subalgebra $\bb$ the subalgebra of $so(d)$ contained in the
subalgebra of upper triangular matrices.  Furthermore we take as a
Cartan subalgebra of $\bb$ the subalgebra $\h$ spanned by the
elements $\tilde{E}_{ii}=E_{ii}-E_{d+1-i,d+1-i}$, for
$i=1,\cdots,k$.  Now a finite-dimensional irreducible
representation of $so(d)$ is determined by its highest weight
$\la\in\h^*$ subject to
\begin{align*}
&\la(\tilde{E}_{ii}-\tilde{E}_{i+1,i+1})\in\Z_+,\allowdisplaybreaks\\
&\la(\tilde{E}_{k-1,k-1}+\tilde{E}_{kk})\in\Z_+,
\end{align*}
for $i=1,\cdots,k-1$.  Let $\la_i=\la(\tilde{E_{ii}})$ and
identify $\la$ with the sequence of complex numbers
$(\la_1,\cdots,\la_k)$. An irreducible representation of $so(2k)$
is finite-dimensional if and only if its highest weight $\la$
satisfies the conditions $\la_1\ge\la_2\cdots\ge\la_k$ with either
$\la_i\in\Z$ or else $\la_i\in\hf+\Z$, $i=1,\cdots,k$ and
$\la_{j}\ge 0$, $j=1,\cdots,k-1$. Furthermore such a weight lifts
to a representation of $SO(d)$ if and only if $\la_i\in\Z_+$.

Let $V$ be a finite-dimensional irreducible $O(d)$-module.  When
regarded as an $so(d)$-module we have the following possibilities:
\begin{itemize}
\item[(i)] $V$ is a direct sum of two irreducible $so(d)$-modules of
highest weights $(\la_1,\la_2,\cdots,\la_k)$ and
$(\la_1,\la_2,\cdots,\la_{k-1},-\la_k)$, respectively, where
$\la_k>0$.
\item[(ii)]
$V$ is an irreducible $so(d)$-module of highest weight
$(\la_1,\la_2,\cdots,\la_{k-1},0)$.
\end{itemize}
Here $\la_i\in\Z_+$ for all $i$. In the first case, that is when
$V$ is the direct sum of the two irreducible $so(d)$-modules we
denote $V$ by $V_{O(d)}^\la$, where we let
$\la=(\la_1,\la_2,\cdots,\la_{k-1},\la_k>0)$.  In the second case
there are two possible choices of $V$, which we denote by
$V_{O(d)}^\la$ and $V_{O(d)}^\la\otimes {\rm det}$, respectively.
Recalling that $O(d)$ is a semidirect product of $SO(d)$ and
$\Z_2$ these two $O(d)$-modules as $SO(d)$-modules are isomorphic.
However as $O(d)$-modules they differ by the determinant
representation so that we may distinguish these two modules as
follows: consider the element $\tau\in O(d)-SO(d)$ that switches
the basis vector $e^{k}$ with $e^{k+1}$ and leaves all other basis
vectors of $\C^d$ invariant. We declare $V^\la_{O(d)}$ to be the
$O(d)$-module on which $\tau$ transforms an $SO(d)$-highest weight
vector trivially.  Note that $\tau$ transforms an $SO(d)$-highest
weight vector in the $O(d)$-module $V^\la_{O(d)}\otimes{\rm det}$
by $-1$.

We may associate Young diagrams to these $O(d)$-highest weights as
follows (cf.~\cite{H2}).  For $\la_1\ge\la_2\cdots\ge\la_k>0$ we
have an obvious Young diagram of length $k$.  When $\la_k=0$, we
associate to the highest weight of $V^\la_{O(d)}$ the usual Young
diagram of length less than $k$. To the highest weight of
$V^\la_{O(d)}\otimes {\rm det}$ we associate the Young diagram
obtained from the Young diagram of $\la$ by replacing its first
column by a column of length $d-\la'_1$.  Here and further, for a
partition $\la$, we denote by $\la'$ its conjugate partition. We
have thus associated to each finite-dimensional irreducible
representation of $O(d)$ a Young diagram $\la$ with
$\la_1'+\la'_2\le d$.

Next consider the case when $d=2k+1$ is odd. We take as a Borel
subalgebra $\bb$ the subalgebra of $so(d)$ spanned by upper
triangular matrices so that a Cartan subalgebra $\h$ of $\bb$ is
again spanned by the elements
$\tilde{E}_{ii}=E_{ii}-E_{d+1-i,d+1-i}$, for $i=1,\cdots,k$.  Now
a finite-dimensional irreducible representation of $so(d)$ is
determined by its highest weight $\la\in\h^*$ subject to
\begin{align*}
&\la(\tilde{E}_{ii}-\tilde{E}_{i+1,i+1})\in\Z_+,\allowdisplaybreaks\\
&\la(\tilde{E}_{kk})\in\hf\Z_+,
\end{align*}
for $i=1,\cdots,k-1$.  We set $\la_i=\la(\tilde{E}_{ii})$ and
identify $\la$ with the sequence of complex numbers
$(\la_1,\la_2,\cdots,\la_k)$.  It follows that a highest weight
$\la$ of $so(2k+1)$ gives a finite-dimensional irreducible
representation if and only $\la_1\ge\la_2\cdots\ge\la_k$ and
$\la_i\in\Z_+$ or else $\la_i\in\hf+\Z_+$, for $i=1,\cdots,k$.

Recall that when $d$ is odd $O(d)$ is a direct product of $SO(d)$
and $\Z_2$. Thus any finite-dimensional irreducible representation
of $O(d)$, when regarded as an $SO(d)$-module, remains
irreducible. Conversely an irreducible representation of $SO(d)$
gives rise to two non-isomorphic $O(d)$-modules that differ from
each other by the determinant representation ${\rm det}$. We let
$V^\la_{O(d)}$ stand for the irreducible $O(d)$-module
corresponding to $\la=(\la_1\ge\la_2\ge\cdots\ge\la_k\ge 0)$ on
which the element $-I$ transforms trivially, so that
$\{V^\la_{O(d)},V^\la_{O(d)}\otimes{\rm det}\}$ with $\la$ ranging
over all partitions as above is a complete set of
finite-dimensional non-isomorphic irreducible $O(d)$-modules.

Similarly as before we may associate Young diagrams to these
$O(d)$-highest weights.  For the highest weight
$\la=(\la_1\ge\la_2\cdots\ge\la_k\ge 0)$ of $V^\la_{O(d)}$ we have
an obvious Young diagram with $l(\la):=\la'_1\le k$. To the
highest weight of $V^\la_{O(d)}\otimes{\rm det}$ we associate the
Young diagram obtained from the Young diagram of $\la$ by
replacing its first column by a column of length $d-\la'_1$.

Let $\epsilon_i\in\h^*$ so that
$\epsilon_i(\tilde{E}_{jj})=\delta_{ij}$.  We put
$x_i=e^{\epsilon_i}$ when dealing with characters of $O(d)$.

\subsection{Irreducible representations of the symplectic group}
Let $d=2k$ and consider the non-degenerate skew-symmetric bilnear
form $<\cdot|\cdot>$ given by the $d\times d$ matrix
\begin{equation*}
\begin{pmatrix}
0&J_{k}\\
-J_{k}&0
\end{pmatrix}.
\end{equation*}
The symplectic group $Sp(d)$ is the subgroup of $GL(d)$ preserving
$<\cdot|\cdot>$. We take as a Borel subalgebra $\bb$ the
subalgebra of $sp(d)$ that is contained in the subalgebra of upper
triangular matrices and a Cartan subalgebra of $\bb$ as the
subalgebra $\h$ spanned by the elements
$\tilde{E}_{ii}=E_{ii}-E_{d+1-i,d+1-i}$, for $i=1,\cdots,k$. A
finite-dimensional irreducible representation of $sp(d)$ is
determined by its highest weight $\la\in\h^*$ subject to
\begin{align*}
&\la(\tilde{E}_{ii}-\tilde{E}_{i+1,i+1})\in\Z_+,\allowdisplaybreaks\\
&\la(\tilde{E}_{k,k})\in\Z_+,
\end{align*}
for $i=1,\cdots,k-1$. As before we let $\la_i=\la(\tilde{E}_{ii})$
and identify $\la$ with the sequence $(\la_1,\la_2,\cdots,\la_k)$.
A highest weight $\la$ of $sp(2k)$ gives a finite-dimensional
irreducible representation if and only if
$\la_1\ge\la_2\cdots\ge\la_k$ and $\la_i\in\Z_+$ for
$i=1,\cdots,k$. Furthermore each such representation lifts to a
unique irreducible representation of $Sp(d)$ and so we obtain an
obvious parameterization of $Sp(d)$-highest weight in terms of
Young diagrams $\la$ with $l(\la)\le \frac{d}{2}$.

We let $\epsilon_i\in\h^*$ so that
$\epsilon_i(\tilde{E}_{jj})=\delta_{ij}$. We put
$y_i=e^{\epsilon_i}$ when dealing with characters of $Sp(2k)$.

\subsection{Irreducible representations of the ortho-symplectic Lie
superalgebra}\label{ospirrep} Let $\C^{m|n}$ be the
$(m|n)$-dimensional complex superspace.  Suppose that $n$ is even
and $(\cdot|\cdot)$ is a supersymmetric non-degenerate bilinear
form, i.e.~it is symmetric on the even subspace $\C^{m|0}$ and
symplectic on the odd subspace $\C^{0|n}$. The orthosymplectic Lie
superalgebra $osp(m|n)$ (cf.~\cite{K}) is defined to be the
subalgebra of $gl(m|n)=gl(m|n)_{\bar{0}}\oplus gl(m|n)_{\bar{1}}$
consisting of those linear transformations preserving the form
$(\cdot|\cdot)$, i.e. $osp(m|n)=osp(m|n)_{\bar{0}}\oplus
osp(m|n)_{\bar{1}}$ with
\begin{equation*}
osp(m|n)_\epsilon=\{A\in
gl(m|n)_{\epsilon}|(Av|w)+(-1)^{\epsilon{\rm deg }v}(v|Aw)=0\},
\end{equation*}
where $v$ and $w$ are any homogeneous vectors of $\C^{m|n}$, ${\rm
deg}v$ here and further denotes the degree of the homogeneous
element $v$ and $\epsilon\in\Z_2$. We will fix the bilinear form
associated to matrix
\begin{equation*}
\begin{pmatrix}
J_m&0&0\\
0&0&J_{n/2}\\
0&-J_{n/2}&0
\end{pmatrix}.
\end{equation*}

We note that $osp(m|n)_{\bar{0}}\cong so(m)\oplus sp(n)$. Let
$\bb$ be a Borel subalgebra of $osp(m|n)$ containing the Borel
subalgebras of $so(m)$ and $sp(n)$ as chosen above so that a
Cartan subalgebra $\h$ of $osp(m|n)$ can be taken to be the
subalgebra spanned by the diagonal matrices
$\tilde{E}_{ii}=E_{ii}-E_{m+1-i,m+1-i}$,
$i=1,\cdots,[\frac{m}{2}]$,
$\tilde{E}_{[\frac{m}{2}]+j,[\frac{m}{2}]+j}=E_{m+j,m+j}-E_{m+n+1-j,m+n+1-j}$,
$j=1,\cdots,\frac{n}{2}$. Here and further the symbol $[r]$ stands
for the largest integer smaller than or equal to $r$. As usual,
highest weight irreducible representations of $osp(m,n)$ are
parameterized by $\la\in\h^*$ and we denote by $\la_i$ the $i$-th
label $\la(\tilde{E}_{ii})$, for
$i=1,\cdots,[\frac{m}{2}]+\frac{n}{2}$. As usual, we will identify
$\la$ with $(\la_1,\la_2,\cdots)$.

Suppose that $m$ is an even integer and consider the following
$\Z$-gradation of $osp(m|n)$. Let $\C^{m|0}=V\oplus V^*$ be a sum
of two isotropic subspaces of $\C^{m|0}$ with respect to the
restriction of the form $(\cdot|\cdot)$ on $\C^{m|0}$. Likewise
let $\C^{0|n}=W\oplus W^*$ be such an isotropic decomposition of
$\C^{0|n}$.  We have $osp(m|n)_{\bar{0}}\cong
S^2(\C^{0|n})\oplus\Lambda^2(\C^{m|0})$ and
$osp(m|n)_{\bar{1}}\cong \C^{m|0}\otimes\C^{0|n}$.  Set
$\G_0=(V\oplus W)\otimes (V\oplus W)^*$,
$\G_1=S^2(V)\oplus\Lambda^2(W)\oplus (V\otimes W)$ and
$\G_{-1}=S^2(V^*)\oplus\Lambda^2(W^*)\oplus (V^*\otimes W^*)$.
This equips $osp(m|n)$ with a $\Z$-gradation with $\G_0$
isomorphic to $gl(\frac{m}{2}|\frac{n}{2})$ such that its standard
Cartan subalgebra is also $\h$.

Now take a finite-dimensional irreducible $\G_0$-module
$V_{\frac{m}{2}|\frac{n}{2}}^{\la}$ of highest weight
$\la\in\h^*$, which we again will identify with a sequence
$(\la_1,\la_2,\cdots)$. We may extend
$V^{\la}_{\frac{m}{2}|\frac{n}{2}}$ trivially to a module over the
parabolic subalgebra $\G_{0}\oplus\G_1$. Inducing it to an
$osp(m|n)$-module, it is clear that it has a unique irreducible
quotient, which we will denote by $V^\la_{osp(m|n)}$. Of course
$V^\la_{osp(m|n)}$ is not finite-dimensional in general. As such
$osp(m|n)$-modules play an important role in the sequel, we will
give a more detailed description of their parameterizations. Let
$\epsilon_i\in\h^*$, $i=1,\cdots,[\frac{m}{2}]+\frac{n}{2}$, be
defined by $\epsilon_i(\tilde{E}_{jj})=\delta_{ij}$. We will label
the simple roots and coroots of $osp(m|n)$ according to the
following diagram.

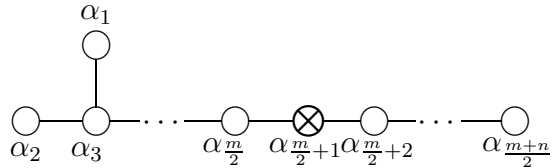
\begin{table}[ha]
\vspace*{-8ex}$
\begin{array}{c c}
\setlength{\unitlength}{0.16in}
\begin{picture}(20,8)
%
%
\put(2.3,4.5){\makebox(0,0)[c]{$\alpha_1$}}
\put(2.3,3.5){\makebox(0,0)[c]{$\bigcirc$}}

\put(2.28,1.45){\line(0,1){1.57}}
\put(0,1){\makebox(0,0)[c]{$\bigcirc$}}
\put(2.3,1){\makebox(0,0)[c]{$\bigcirc$}}
\put(6.85,1){\makebox(0,0)[c]{$\bigcirc$}}
\put(9.25,1){\makebox(0,0)[c]{$\bigotimes$}}
\put(11.4,1){\makebox(0,0)[c]{$\bigcirc$}}
\put(16,1){\makebox(0,0)[c]{$\bigcirc$}}
\put(0.45,1){\line(1,0){1.4}} \put(2.72,1){\line(1,0){1}}
\put(5.2,1){\line(1,0){1.2}} \put(7.28,1){\line(1,0){1.45}}
\put(9.7,1){\line(1,0){1.25}} \put(11.85,1){\line(1,0){0.9}}
\put(14.25,1){\line(1,0){1.3}}
\put(4.5,0.95){\makebox(0,0)[c]{$\cdots$}}
\put(13.5,0.95){\makebox(0,0)[c]{$\cdots$}}
\put(0,0){\makebox(0,0)[c]{$\alpha_2$}}
\put(2,0){\makebox(0,0)[c]{$\alpha_3$}}
\put(6.5,0){\makebox(0,0)[c]{$\alpha_{{\frac{m}{2}}}$}}
\put(9.15,0){\makebox(0,0)[c]{$\alpha_{\frac{m}{2}+1}$}}
\put(11.5,0){\makebox(0,0)[c]{$\alpha_{\frac{m}{2}+2}$}}
\put(16.1,0){\makebox(0,0)[c]{$\alpha_{\frac{m+n}{2}}$}}
\end{picture}
\end{array}$
\end{table}

\noindent Here $\alpha_1=-\epsilon_1-\epsilon_2,
\alpha_2=\epsilon_1-\epsilon_2,\cdots,
\alpha_{\frac{m}{2}}=\epsilon_{\frac{m}{2}-1}-\epsilon_{\frac{m}{2}},
\cdots,\alpha_{\frac{m+n}{2}}=\epsilon_{\frac{m+n}{2}-1}-
\epsilon_{\frac{m+n}{2}}$, and, as is customary, $\bigotimes$
denotes an isotropic root. Thus if
$\la=(\la_1,\la_2,\cdots,\la_{\frac{m+n}{2}})$ is the highest
weight of a finite-dimensional irreducible
$gl(\frac{m}{2}|\frac{n}{2})$-module
$V^\la_{\frac{m}{2}|\frac{n}{2}}$, then the labels of the
irreducible highest weight module $V^\la_{osp(m|n)}$ with respect
to the above Dynkin diagram is given by
\begin{equation}\label{osplabels}
(-\la_1-\la_2,\la_1-\la_2,\cdots,\la_{\frac{m}{2}-1}-\la_{\frac{m}{2}},
\la_{\frac{m}{2}}+\la_{\frac{m}{2}+1},\la_{\frac{m}{2}+1}-\la_{\frac{m}{2}+2},\cdots).
\end{equation}
 When dealing with characters of
$osp(m|n)$ we will use the notation $x_j=e^{\epsilon_j}$, for
$j=1,\cdots,\frac{m}{2}$ and $z_l=e^{\epsilon_{\frac{m}{2}+l}}$,
for $l=1,\cdots,\frac{n}{2}$.

On the superspace $\C^{m|n}$ with $m$ even we may take a
skew-supersymmetric non-degenerate bilinear form $(\cdot|\cdot)$,
i.e.~it is symplectic on the even subspace $\C^{m|0}$ and
symmetric on the odd subspace $\C^{0|n}$. In the same fashion we
may define the symplectic-orthogonal Lie superalgebra $spo(m|n)$
to be the subalgebra of $gl(m|n)$ preserving $(\cdot|\cdot)$.  We
remark that as Lie superalgebras we have $spo(m|n)\cong osp(n|m)$
and hence our discussion of the ortho-symplectic Lie superalgebra
carries over to $spo(m|n)$, for $n$ even, with minor modification.
We label the simple roots and coroots according to the following
diagram.

\begin{table}[hb]
\vspace*{-8ex}$
\begin{array}{c c}
\setlength{\unitlength}{0.16in}
\begin{picture}(20,5)
%
%
%
\put(0,1){\makebox(0,0)[c]{$\bigcirc$}}
\put(2.4,1){\makebox(0,0)[c]{$\bigcirc$}}
\put(6.85,1){\makebox(0,0)[c]{$\bigcirc$}}
\put(9.25,1){\makebox(0,0)[c]{$\bigotimes$}}
\put(11.4,1){\makebox(0,0)[c]{$\bigcirc$}}
\put(16,1){\makebox(0,0)[c]{$\bigcirc$}}
\put(0.35,0.75){$\Longrightarrow$} \put(2.82,1){\line(1,0){0.8}}
\put(5.2,1){\line(1,0){1.2}} \put(7.28,1){\line(1,0){1.45}}
\put(9.7,1){\line(1,0){1.25}} \put(11.81,1){\line(1,0){0.9}}
\put(14.25,1){\line(1,0){1.28}}
\put(4.5,0.95){\makebox(0,0)[c]{$\cdots$}}
\put(13.5,0.95){\makebox(0,0)[c]{$\cdots$}}
\put(0,0){\makebox(0,0)[c]{$\alpha_1$}}
\put(2.4,0){\makebox(0,0)[c]{$\alpha_2$}}
\put(6.5,0){\makebox(0,0)[c]{$\alpha_{\frac{m}{2}}$}}
\put(9.15,0){\makebox(0,0)[c]{$\alpha_{\frac{m}{2}+1}$}}
\put(11.5,0){\makebox(0,0)[c]{$\alpha_{\frac{m}{2}+2}$}}
\put(16.1,0){\makebox(0,0)[c]{$\alpha_{\frac{m+n}{2}}$}}
\end{picture}
\end{array}$
\end{table}

\noindent Here $\alpha_1=-2\epsilon_1,
\alpha_2=\epsilon_1-\epsilon_2,\cdots,
\alpha_{\frac{m}{2}}=\epsilon_{\frac{m}{2}-1}-\epsilon_{\frac{m}{2}},
\cdots,\alpha_{\frac{m+n}{2}}=\epsilon_{\frac{m+n}{2}-1}-
\epsilon_{\frac{m+n}{2}}$.  Similarly we will denote the
irreducible quotient of the induced $gl(m|n)$-module
$V^\la_{\frac{m}{2}|\frac{n}{2}}$ by $V^\la_{spo(m|n)}$. So if
$\la=(\la_1,\la_2,\cdots,\la_{\frac{m+n}{2}})$ is the
$gl(\frac{m}{2}|\frac{n}{2})$-labels of
$V^\la_{\frac{m}{2}|\frac{n}{2}}$, then the $spo(m|n)$-labels of
$V^\la_{spo(m|n)}$ are
\begin{equation}\label{spolabels}
(-\la_1,\la_1-\la_2,\cdots,\la_{\frac{m}{2}-1}-\la_{\frac{m}{2}},
\la_{\frac{m}{2}}+\la_{\frac{m}{2}+1},\la_{\frac{m}{2}+1}-\la_{\frac{m}{2}+2},\cdots).
\end{equation}
When dealing with characters of $spo(m|n)$ we will use the
notation $y_j=e^{\epsilon_j}$, for $j=1,\cdots,\frac{m}{2}$ and
$z_l=e^{\epsilon_{\frac{m}{2}+l}}$, for $l=1,\cdots,\frac{n}{2}$.

\section{The $(gl(d),gl(m|n))$-duality}\label{glglduality}

In this section we present some results on
$(gl(d),gl(m|n))$-duality that will be used later on. In
particular, \thmref{basis} constructs explicit bases for
irreducible $gl(d)\times gl(m|n)$-modules appearing in the
decomposition of $S(\C^d\otimes\C^{m|n})$, and we believe the
result to be new.

Consider the natural actions of $gl(d|q)$ on $\C^{d|q}$ and
$gl(m|n)$ on $\C^{m|n}$.  We can form the $gl(d|q)\times
gl(m|n)$-module $\C^{d|q}\otimes\C^{m|n}$.  We have an induced
action on the symmetric tensor $S(\C^{d|q}\otimes\C^{m|n})$. This
action is completely reducible and in fact $(gl(d|q),gl(m|n))$ is
a dual pair in the sense of Howe \cite{CW1} (see also \cite{S1}).
Since in this paper we will only concern ourselves with the case
when $q=0$, we will make this assumption in what follows. In this
case we have the following decomposition
\begin{equation}\label{glgl-duality}
S(\C^d\otimes\C^{m|n})\cong\sum_{\la}V^\la_d\otimes V^\la_{m|n},
\end{equation}
The sum in \eqref{glgl-duality} is over all partitions of integers
$\la=(\la_1,\la_2,\cdots,\la_d)$ of length $l(\la)$ not exceeding
$d$ subject to $\la_{m+1}\le n$.  Since $l(\la)\le d$ we may
regard $\la$ as a highest weight for an irreducible $gl(d)$-module
so that there is no ambiguity in $V^\la_d$.  The meaning of
$V^\la_{m|n}$ as a $gl(m|n)$-module was explained in
\secref{irrepglmn}.

In the sequel it is important to have an explicit formula for the
joint highest weight vectors of the irreducible component
$V^\la_d\otimes V^\la_{m|n}$ in \eqnref{glgl-duality}. (See also
\cite{OP} and \cite{N} for different descriptions of these
vectors.) In order to present them we need to introduce some more
notation.

We let $e^1,\ldots,e^d$ denote the standard basis for the standard
$gl(d)$-module.  Similarly we let $e_1,\ldots,e_m;f_1,\ldots,f_n$
denote the standard homogeneous basis for the standard
$gl(m|n)$-module. The weights of $e^i$, $e_l$ and $f_k$ are
denoted by $\tilde{\epsilon}_i$, $\epsilon_l$ and $\delta_k$, for
$1\le i\le d$, $1\le l\le m$ and $1\le k\le n$, respectively. We
set
\begin{equation}\label{generators}
x_l^i:=e^i\otimes e_l,\quad \eta_k^i:=e^i\otimes f_k.
\end{equation}

We will denote by $\C[\x,\etabf]$ the polynomial superalgebra
generated by \eqnref{generators}. By identifying
$S(\C^{d}\otimes\C^{m|n})$ with the polynomial superalgebra
$\C[\x,\etabf]$ the commuting pair $(gl(d),gl(m|n))$ may be
realized as first order differential operators as follows ($1\le
i,i'\le d$, $1\le s,s'\le m$ and $1\le k,k'\le n$):
\begin{align}
&\sum_{j=1}^m x_{j}^{i}\frac{\partial}{\partial
x_j^{i'}}+\sum_{j=1}^n\eta_j^{i}\frac{\partial}{\partial\eta_j^{i'}},\label{glpq} \allowdisplaybreaks\\
&\sum_{j=1}^d x_{s}^{j}\frac{\partial}{\partial x_{s'}^j},\quad
\sum_{j=1}^d\eta_{k'}^{j}\frac{\partial}{\partial\eta_k^j},\quad
\sum_{j=1}^d x_s^j\frac{\partial}{\partial\eta_k^j},\quad
\sum_{j=1}^d\eta_k^j\frac{\partial}{\partial x_s^j}.\label{glmn1}
\end{align}
Here \eqnref{glpq} spans a copy of $gl(d)$, while \eqnref{glmn1}
spans a copy of $gl(m|n)$.

The standard Cartan subalgebras of $gl(d)$ and $gl(m|n)$ are
spanned, respectively, by
\begin{equation*}
\sum_{j=1}^m x_{j}^{i}\frac{\partial}{\partial
x_j^{i}}+\sum_{j=1}^n\eta_j^{i}\frac{\partial}{\partial\eta_j^{i}}
\quad {\rm and} \quad \sum_{j=1}^d
x_{s}^{j}\frac{\partial}{\partial x_{s}^j},\
\sum_{j=1}^d\eta_{k}^{j}\frac{\partial}{\partial\eta_k^j},
\end{equation*}
while the nilpotent radicals are respectively generated by the
simple root vectors

\begin{eqnarray*}
\sum_{j=1}^m x_{j}^{i-1}\frac{\partial}{\partial
x_j^{i}}+\sum_{j=1}^n\eta_j^{i-1}\frac{\partial}{\partial\eta_j^{i}},
\quad 1< i \le d,
\end{eqnarray*}
and
\begin{align*}
\sum_{j=1}^d x_{s-1}^{j}\frac{\partial}{\partial x_{s}^j},\
\sum_{j=1}^d\eta_{k-1}^{j}\frac{\partial}{\partial\eta_k^j},\
 \sum_{j=1}^d x_m^j\frac{\partial}{\partial\eta_1^j},\quad 1< s \le m,1< k \le
n.
\end{align*}

We will consider two separate cases, namely $m\ge d$ and $m<d$.

First suppose that $m\ge d$.  Here the condition $\la_{m+1}\le n$
is vacuous. For $1\le r\le {\rm min}(d,m)$ define
\begin{equation}\label{deltar}
\Delta_r:={\rm det}\begin{pmatrix}
x_{1}^{1}&x_{2}^{1}&\cdots&x_{r}^{1}\\
x_{1}^{2}&x_{2}^{2}&\cdots&x_{r}^{2}\\
\vdots&\vdots&\vdots&\vdots\\
x_{1}^{r}&x_{2}^{r}&\cdots&x_{r}^{r}\\
\end{pmatrix}.
\end{equation}

\begin{thm}\cite{CW1}\label{aux2}
In the case when $m\ge d$, the $gl(d) \times gl(m|n)$ highest
weight vectors in $ \C[\x, \etabf]$ associated to the weight $\la$
is given by the product
$\Delta_{\la'_1}\Delta_{\la'_2}\cdots\Delta_{\la'_{\la_1}}$.
\end{thm}

We now consider the case $d>m$.  It is readily checked that the
highest weight vectors associated to Young diagrams $\la$ with
$\la_{m+1}=0$ can be obtained just as in the previous case so that
we may assume that $l(\la)>m$. Let
$\lambda_1',\lambda_2',\ldots,\la'_{\la_1}$ denote its column
lengths as usual.  We have $d \ge \la'_1 \ge \la'_2 \ldots \ge
\la'_{\la_1}$ and $m\ge\la'_{n+1}$.  For $m< r\le d$, consider the
following determinant of an $r\times r$ matrix:
\begin{equation} \label{eq_det}
\Delta_{k,r}:={\rm det}
\begin{pmatrix}
x_1^1&x_1^2&\cdots &x_1^r\\
x_2^1&x_2^2&\cdots &x_2^r\\
\vdots&\vdots&\cdots &\vdots\\
x_m^1&x_m^2&\cdots &x_m^r\\
\eta_k^1&\eta_k^2&\cdots &\eta_k^r\\
\eta_k^1&\eta_k^2&\cdots &\eta_k^r\\
\vdots&\vdots&\cdots &\vdots\\
\eta_k^1&\eta_k^2&\cdots &\eta_k^r\\
\end{pmatrix},\quad k=1,\ldots,n.
\end{equation}
That is, the first $m$ rows are filled by the vectors $(x_j^1,
\ldots, x_j^r)$, for $j=1, \ldots, m$, in increasing order and the
last $r-m$ rows are filled with the same vector
$(\eta_k^1,\ldots,\eta_k^r)$. Here the determinant of a matrix
$$
 A :=\begin{pmatrix} a_1^1 & a_1^2 &\cdots&a_1^r\\
a_2^1&a_2^2&\cdots&a_2^r\\
\vdots&\vdots&\cdots &\vdots\\
a_r^1&a_r^2&\cdots&a_r^r\\
\end{pmatrix},
$$
with matrix entries possibly involving Grassmann variables
$\eta_k^i$, is by definition the expression $\sum_{\sigma\in
S_r}(-1)^{l(\sigma)}a_1^{\sigma(1) }a_2^{\sigma (2) }\cdots
a_r^{\sigma (r) }$, where $l(\sigma)$ is the length of $\sigma$ in
the symmetric group $S_r$.

\begin{thm}\label{glpmn-duality}\cite{CW1}
In the case when $m< d$, the $gl(d) \times gl(m|n)$ highest weight
vectors in $ \C[\x, \etabf]$ associated to the weight $\la$ is
given by the product
\begin{equation}\label{hwvform}
\prod_{k=1}^\nu\Delta_{k,\la'_k}\prod_{j=r+1}^{\la_1}\Delta_{\la'_j},
\end{equation}
where $\nu$ is defined by $\la'_\nu>m$ and $\la'_{\nu+1}\le m$.
\end{thm}

For application purposes it is useful to construct an explicit
basis for the $gl(d)\times gl(m|n)$-modules $V_{m|n}^\la$ that
appear in the decomposition of $S(\C^d\otimes\C^{m|n})$.  This we
will do now.

Recall that $\la$ is a partition (or a Young diagram) which lies
in the $(m|n)$-hook of length not exceeding $d$. Let
$x_1,\cdots,x_m$ and $\eta_1,\cdots,\eta_n$ be even and odd
indeterminates, respectively. We form a tableau of shape $\la$ by
filling the nodes of $\la$ from the set
$\{x_{1},\cdots,x_{m},\eta_1,\cdots,\eta_n\}$ so that the
resulting tableau $T$ is $(m|n)$-semi-standard.  This means that
we first fill the nodes of a sub-diagram $\mu\subseteq\la$ with
the even indeterminates $\{x_1,\cdots,x_m\}$ so that the resulting
sub-tableau is semi-standard.  Then we fill the skew-diagram
$\la/\mu$ with odd indeterminates $\{\eta_1,\cdots,\eta_n\}$ so
that its transpose is semi-standard. Let us suppose that the
$i$-th column of $T$ has length $r$ and is filled from top to
bottom by
\begin{equation}\label{ithcol}
(x_{i_1},x_{i_2},\cdots,x_{i_s},\eta_{j_1},\cdots,\eta_{j_t}).
\end{equation}
We associate to \eqnref{ithcol} the following determinant
\begin{equation} \label{ithcoldet}
\Delta_i^T:={\rm det}
\begin{pmatrix}
x_{i_1}^1&x_{i_1}^2&\cdots &x_{i_1}^r\\
x_{i_2}^1&x_{i_2}^2&\cdots &x_{i_2}^r\\
\vdots&\vdots&\cdots &\vdots\\
x_{i_s}^1&x_{i_s}^2&\cdots &x_{i_s}^r\\
\eta_{j_1}^1&\eta_{j_1}^2&\cdots &\eta_{j_1}^r\\
\eta_{j_2}^1&\eta_{j_2}^2&\cdots &\eta_{j_2}^r\\
\vdots&\vdots&\cdots &\vdots\\
\eta_{j_t}^1&\eta_{j_t}^2&\cdots &\eta_{j_t}^r\\
\end{pmatrix},
\end{equation}
where $r=s+t$. We set $\Delta^T=\prod_{i=1}^{\la_1}\Delta_i^T$.

\begin{thm}\label{glhwv}
The set $\{\Delta^T\}$, with $T$ running over all
$(m|n)$-semi-standard tableaux of shape $\la$, is a basis for the
space of $gl(d)$-highest weight vectors in
$S(\C^d\otimes\C^{m|n})$ of highest weight $\la$.
\end{thm}

\begin{proof}
It is easy to see that every $\Delta^T$ is a $gl(d)$-highest
weight vector of $gl(d)$-highest weight $\la$. Now according to
\cite{BR} the dimension of $V^\la_{m|n}$ equals the number of
$(m|n)$-semi-standard tableaux of shape $\la$ and hence it is
enough to show that the set $\{\Delta^T\}$ is a linearly
independent set. Now due to weight considerations it is enough to
prove that the set of $\{\Delta^T\}$, where $T$ is over all
$(m|n)$-semi-standard tableaux with fixed occurrence of
$\{x_1,\cdots,x_m,\eta_1,\cdots,\eta_n\}$, is linearly
independent.  We proceed by induction on the number of odd
indeterminates that occur inside the $T$'s. If that number is
zero, then the conclusion of the theorem is know to be true (see
e.g.~\cite{FH}). Thus we may assume that at least one odd
indeterminate occurs in all of the $T$'s.

Now let $\eta_i$ be the odd indeterminate appearing in all
$\Delta^T$ with $i$ minimal.  Let
\begin{equation}\label{sumla}
\sum_{T}\la_T\Delta^T=0.
\end{equation}
We embed $S(\C^{d}\otimes\C^{m|n})$ into
$S(\C^d\otimes\C^{m+1|n})$ so that we may regard \eqnref{sumla} as
a sum in $S(\C^d\otimes\C^{m+1|n})$.  We apply to \eqnref{sumla}
the linear map
\begin{equation*}
A=\sum_{j=1}^d x^j_{m+1}\frac{\pa}{\pa\eta_i^j}.
\end{equation*}
It is clear that the resulting sum is of the form
\begin{equation*}
\sum_{T}\sum_{S\in\La_T}\la_S\Delta^S,
\end{equation*}
where $\La_T$ is the set of all tableaux obtained from $T$ by
replacing one of the $\eta_i$-nodes by an $x_{m+1}$-node.  We may
assume that all $S$ are $(m+1|n)$-semi-standard with one less odd
node. Furthermore each $\la_S$ is a non-zero positive integral
multiple of $\la_T$.  (Note that $\la_S=p\la_S$ if and only if
$\eta_i$ appears with multiplicity $p$ in some column and $S$ is
obtained from $T$ by replacing the first $\eta_i$ node of this
column by $x_{m+1}$.)

We claim that all $S$ are distinct $(m+1|n)$-semi-standard
tableaux and thus by induction they are linearly independent.
This implies $\la_S=0$ and hence $\la_T=0$ and we are done.

In order to prove the claim we consider two cases.

In the first case suppose that $S$ and $S'$ are obtained from the
same $T$.  But in this case $S$ and $S'$ are obviously different,
since $S$ and $S'$ are obtained from $T$ by replacing $\eta_i$ by
$x_{m+1}$ in different columns.

Now suppose that $S$ and $S'$ are obtained from $T$ and $T'$,
respectively, and $T\not=T'$.  If the positions of $\eta_i$ in $T$
and $T'$ are the same, then $T$ and $T'$ differ at some $\eta_s$
node, $i\not=s$.  But then $S$ and $S'$ also differ at this
particular $\eta_s$-node as well.  If on the other hand $T$ and
$T'$ differ at some $\eta_i$ node, then this means that $T$ at a
node has $\eta_i$, while at the same node $T'$ has some $\eta_s$,
$i\not=s$, or $x_k$, $k\not=m+1$.  But then in all of $S'$ this
particular node is always $\eta_s$ or $x_k$, while in all $S$ this
particular node is either $\eta_i$ or $x_{m+1}$.  Thus $S$ and
$S'$ must be distinct.
\end{proof}

Let $\la$ be a Young diagram lying in the $(m|n)$-hook of with
$l(\la)\le d$ and $T$ be an $(m|n)$-semi-standard tableau of shape
$\la$. We may fill the boxes of the Young diagram $\la$ from the
set $\{x^1,\cdots,x^d\}$ in a way so that the resulting tableau
$T'$ is semi-standard.  Let the $i$-th column of $T'$ be filled by
$\{x^{k_1},\cdots,x^{k_r}\}$. Suppose that a joint $gl(m|n)\times
gl(d)$-highest weight vector is of the form \eqnref{hwvform}.  We
may replace the upper indices $1,2,\cdots,r$ of all the entries in
$\Delta_{i,\la'_i}$ (or $\Delta_{\la'_i}$) by
$k_1,k_2,\cdots,k_r$. Let us call the resulting determinant
$\Delta^{T'}_{i,\la'_i}$ (or $\Delta^{T'}_{\la_i'}$) and consider
the following product of determinant.
\begin{equation*}
\prod_{k=1}^v\Delta^{T'}_{k,\la'_k}\prod_{j=v+1}^{\la_1}\Delta^{T'}_{\la'_j}.
\end{equation*}
It is clear from symmetry between the upper and lower indices that
when $T'$ ranges over all semi-standard tableaux we obtain a basis
for the $gl(m|n)$-highest weight vectors of highest weight $\la$
in $S(\C^d\otimes\C^{m|n})$.

Now let $\la$ be a Young diagram lying in the $(m|n)$-hook of with
$l(\la)\le d$ and $T$ and $T'$ as before.  Let the $i$-th column
of $T'$ and $T$ be filled by
\begin{align*}
&\{x^{k_1},\cdots,x^{k_r}\}\\
&(x_{i_1},x_{i_2},\cdots,x_{i_s},\eta_{j_1},\cdots,\eta_{j_t}),
\end{align*}
respectively, from top to the bottom. To the $i$-th column of the
bi-tableau $(T,T')$ we associate the following determinant
\begin{equation} \label{ithcolbasis}
\Delta_i^{(T,T')}:={\rm det}
\begin{pmatrix}
x_{i_1}^{k_1}&x_{i_1}^{k_2}&\cdots &x_{i_1}^{k_r}\\
x_{i_2}^{k_1}&x_{i_2}^{k_2}&\cdots &x_{i_2}^{k_r}\\
\vdots&\vdots&\cdots &\vdots\\
x_{i_s}^{k_1}&x_{i_s}^{k_2}&\cdots &x_{i_s}^{k_r}\\
\eta_{j_1}^{k_1}&\eta_{j_1}^{k_2}&\cdots &\eta_{j_1}^{k_r}\\
\eta_{j_2}^{k_1}&\eta_{j_2}^{k_2}&\cdots &\eta_{j_2}^{k_r}\\
\vdots&\vdots&\cdots &\vdots\\
\eta_{j_t}^{k_1}&\eta_{j_t}^{k_2}&\cdots &\eta_{j_t}^{k_r}\\
\end{pmatrix},
\end{equation}
where again $r=s+t$. We set
$\Delta^{(T,T')}=\prod_{i=1}^{\la_1}\Delta_i^{(T,T')}$.  The
following theorem gives an explicit basis for each irreducible
$gl(d)\times gl(m|n)$-component in $S(\C^d\otimes\C^{m|n})$.

\begin{thm}\label{basis}
The set $\Delta^{(T,T')}$, where $T'$ is semi-standard in
$\{x^1,\cdots,x^d\}$ and $T$ is $(m|n)$-semi-standard in
$\{x_1,\cdots,x_m,\eta_1,\cdots,\eta_n\}$, is a basis for
$V^\la_d\otimes V^\la_{m|n}$ in $S(\C^d\otimes\C^{m|n})$.
\end{thm}

\begin{proof} Given $\Delta^{(T,T')}$ with $(T,T')$ fixed.
By \thmref{glhwv} and the Jacobson density theorem (more precisely
by Burnside's theorem) we can find an element $(a\otimes b)\in
U(gl(d))\otimes U(gl(m|n))$ such that $(a\otimes
b)\Delta^{(T,T')}$ is the joint $gl(d)\times gl(m|n)$-highest
weight vector and $a\otimes b$ annihilates all $\Delta^{(S,S')}$,
for $(S,S')\not=(T,T')$.  This implies that the set
$\{\Delta^{(T,T')}\}$ is linearly independent. But the number of
semi-standard tableaux in $\{x^1,\cdots,x^d\}$ times the number of
$(m|n)$-semi-standard tableaux in
$\{x_1,\cdots,x_m,\eta_1,\cdots,\eta_n\}$ is precisely the
dimension of the space $V^\la_d\otimes V^\la_{m|n}$.
\end{proof}

\begin{rem}
The above theorem is known in the case when $n=0$ (see
e.g.~\cite{FH}).
\end{rem}

\section{The $(O,spo)$- and $(Sp,osp)$-duality}\label{duality}

Let $\C^d$ be the $d$-dimensional complex vector space with
standard basis $\{e^1,e^2,\allowbreak \cdots,e^d\}$. Let $O(d)$ be
the orthogonal group leaving invariant the symmetric bilinear form
$(\cdot|\cdot)$ as in \secref{parameterization}, and let
$\C^{m|n}$ be the superspace of dimension $(m|n)$. The natural
action of $O(d)$ on $\C^d$ extends to an action on
$\C^d\otimes\C^{m|n}$. This action gives rise to an action of
$O(d)$ on the supersymmetric tensor $S(\C^d\otimes\C^{m|n})$,
which, as in \secref{glglduality}, we identify with
$\C[\x,\etabf]$, the commutative superalgebra in
\eqnref{generators}.  As the action $gl(d)$ under this
identification gets identified with certain first order
differential operators as in \eqnref{glpq}, the action of the Lie
algebra of $O(d)$ naturally gets identified with certain first
order differential operators as well.

Consider the following first order differential operators
\begin{align}
E^{xx}_{is}=\sum_{j=1}^d x^j_i\frac{\partial}{\partial
x_s^j}+\frac{d}{2}\delta_{is},\quad
E^{x\eta}_{ik}=\sum_{j=1}^d x^j_i\frac{\partial}{\partial \eta_k^j}, \label{glmn}\allowdisplaybreaks\\
E^{\eta x}_{ki}=\sum_{j=1}^d \eta^j_k\frac{\partial}{\partial
x_i^j},\quad E^{\eta\eta}_{tk}=\sum_{j=1}^d
\eta^j_t\frac{\partial}{\partial
\eta_k^j}-\frac{d}{2}\delta_{ik},\nonumber
\end{align}
where $i,s=1,\cdots,m$ and $k,t=1,\cdots,n$. It is evident that
they form a basis for the Lie superalgebra $gl(m|n)$ and it is
clear that $O(d)$ commutes with $gl(m|n)$.

Next consider another set of operators on $\C[\x,\etabf]$.
\begin{align*}
&{}^OI^{xx}_{is}=\sum_{j=1}^d x^j_i x^{d+1-j}_s,\quad
{}^OI^{x\eta}_{ik}=\sum_{j=1}^d x^j_i
\eta^{d+1-j}_k,\quad ^{O}I^{\eta\eta}_{kt}=\sum_{j=1}^d \eta^j_k \eta^{d+1-j}_t,
\allowdisplaybreaks\\
&{}^O\Delta^{xx}_{is}=\sum_{j=1}^d \frac{\partial}{\partial x^j_i}
\frac{\partial}{\partial x^{d+1-j}_s},\
{}^O\Delta^{x\eta}_{ik}=\sum_{j=1}^d \frac{\partial}{\partial
x^j_i} \frac{\partial}{\partial \eta^{d+1-j}_k},\
{}^O\Delta^{\eta\eta}_{kt}=\sum_{j=1}^d \frac{\partial}{\partial
\eta^j_k} \frac{\partial}{\partial \eta^{d+1-j}_t},
\end{align*}
where $1\le i\le s\le m$ and $1\le k<t\le n$.  We note that these
operators also commute with the action of $O(d)$ on
$\C[\x,\etabf]$. It is not hard to see that these operators
together with \eqnref{glmn} form a basis of the
symplectic-orthogonal Lie superalgebra $spo(2m|2n)$.  In fact,
using the $\Z$-gradation of $spo(2m|2n)$ given in
\secref{parameterization}, we have
$\G_{1}=\sum\C{}^O\Delta^{xx}_{is}+\sum\C{}^O\Delta^{x\eta}_{ik}
+\sum\C{}^O\Delta^{\eta\eta}_{kl}$ and $\G_{-1}=\sum\C
{}^OI^{xx}_{is}+\sum\C {}^OI
^{x\eta}_{ik}+\sum\C{}^OI^{\eta\eta}_{kl}$. Thus on
$\C[\x,\etabf]$ we have an action of $O(d)\times spo(2m|2n)$.

An element $f\in\C[\x,\etabf]$ will be called {\em
${}^O\Delta$-harmonic}, if
$^O\Delta^{xx}_{is}f=^O\Delta^{x\eta}_{ik}f=
{}^O\Delta^{\eta\eta}_{kl}f=0$. The space of
${}^O\Delta$-harmonics will be denoted by ${}^OH$. Note that since
$[gl(m|n),\G_1]\subseteq\G_1$ the space ${}^OH$ is invariant under
the action of $gl(m|n)$. Also ${}^OH$ is clearly invariant under
the action of $O(d)$.  Hence we have an action of $O(d)\times
gl(m|n)$ on ${}^OH$. Let ${}^OI$ be the subalgebra of
$\C[\x,\etabf]$ generated by ${}^OI^{xx}_{is}$, ${}^OI
^{x\eta}_{ik}$ and ${}^OI^{\eta\eta}_{kl}$. It is clear that
${}^OI$ is the subalgebra of $O(d)$-invariants in $\C[\x,\etabf]$.
We have the following theorem.

\begin{thm}\label{O-duality}\cite{H1}
The pairs $(O(d), spo(2m|2n))$ and $(O(d), gl(m|n))$ form dual
reductive Howe pairs on $S(\C^d\otimes\C^{m|n})$ and on ${}^OH$,
respectively. Thus we have
\begin{align*}
&\C[\x,\etabf]\cong\sum_{\la}V^\la_{O(d)}\otimes V^{\la'}_{spo(2m|2n)},\\
&{}^OH\cong\sum_{\la}V^\la_{O(d)}\otimes V^{\la''}_{m|n},
\end{align*}
where $\la$ is summed over a set of irreducible $O(d)$-highest
weights.  Here $\la'$ and $\la''$ are certain non-isomorphic
irreducible $spo(2m|2n)$- and $gl(m|n)$-highest weights,
respectively. Furthermore the map
${}^OI\otimes{}^OH\rightarrow\C[\x,\etabf]$ given by
multiplication is surjective and we have, for each $\la$,
$V^{\la'}_{spo(2m|2n)}={}^OIV^{\la''}_{m|n}$.
\end{thm}

Let $d$ be an even integer and consider the $d$-dimensional
complex vector space with the standard basis $e^1,e^2,\cdots,e^d$
and equipped with the non-degenerate skew-symmetric bilinear form
$<\cdot|\cdot>$ as in \secref{parameterization}. Let $Sp(d)$ be
the corresponding symplectic group. Again we have an action of
$Sp(d)$ on $\C^d\otimes\C^{m|n}$, inducing an action of $Sp(d)$ on
the supersymmetric tensor $S(\C^d\otimes\C^{m|n})$, which we again
identify with $\C[\x,\etabf]$.

Introduce the following operators
\begin{align*}
&{}^{Sp}I^{xx}_{is}=\sum_{j=1}^{\frac{d}{2}} \Big{(}x^j_i
x^{d+1-j}_s-x^{d+1-j}_i x^j_s\Big{)},\quad
{}^{Sp}I^{x\eta}_{ik}=\sum_{j=1}^{\frac{d}{2}} \Big{(}x^j_i
\eta^j_k-x^{d+1-j}_i
\eta^{d+1-j}_k\Big{)},\allowdisplaybreaks\\
&{}^{Sp}I^{\eta\eta}_{kt}=\sum_{j=1}^{\frac{d}{2}} \Big{(}\eta^j_k
\eta^{d+1-j}_t-\eta^{d+1-j}_k \eta^j_t\Big{)}, \
{}^{Sp}\Delta^{xx}_{is}=\sum_{j=1}^{\frac{d}{2}}
\Big{(}\frac{\partial}{\partial x^j_i} \frac{\partial}{\partial
x^{d+1-j}_s}-\frac{\partial}{\partial x^{d+1-j}_i}
\frac{\partial}{\partial
x^j_s}\Big{)},\allowdisplaybreaks\\
&{}^{Sp}\Delta^{x\eta}_{ik}=\sum_{j=1}^{\frac{d}{2}}
\Big{(}\frac{\partial}{\partial x^j_i} \frac{\partial}{\partial
\eta^{d+1-j}_k}-\frac{\partial}{\partial x^{d+1-j}_i}
\frac{\partial}{\partial \eta^j_k}\Big{)},\allowdisplaybreaks\\
&{}^{Sp}\Delta^{\eta\eta}_{kt}=\sum_{j=1}^{\frac{d}{2}}
\Big{(}\frac{\partial}{\partial \eta^j_k} \frac{\partial}{\partial
\eta^{d+1-j}_t}-\frac{\partial}{\partial \eta^{d+1-j}_k}
\frac{\partial}{\partial \eta^j_t}\Big{)},
\end{align*}
where $1\le i< s\le m$ and $1\le k\le t\le n$. It is again not
hard to see that these operators together with \eqnref{glmn} form
a basis for the Lie superalgebra $osp(2m|2n)$ and their actions
and that of $Sp(d)$ on $\C[\x,\etabf]$ commute.

An element $f\in\C[\x,\etabf]$ will be called {\em
${}^{Sp}\Delta$-harmonic}, if
$^{Sp}\Delta^{xx}_{is}f=^{Sp}\Delta^{x\eta}_{ik}f=
{}^{Sp}\Delta^{\eta\eta}_{kl}f=0$. The space of
${}^{Sp}\Delta$-harmonics will be denoted by ${}^{Sp}H$. Similarly
we have an action of $Sp(d)\times gl(m|n)$ on ${}^{Sp}H$. Let
${}^{Sp}I$ be the subalgebra of $\C[\x,\etabf]$ generated by
${}^{Sp}I^{xx}_{is}$, ${}^{Sp}I ^{x\eta}_{ik}$ and
${}^{Sp}I^{\eta\eta}_{kl}$ so that ${}^{Sp}I$ is the subalgebra of
$Sp(d)$-invariants in $\C[\x,\etabf]$. In a similar fashion we
have the following theorem.

\begin{thm}\label{Sp-duality}\cite{H1}
The pairs $(Sp(d), osp(2m|2n))$ and $(Sp(d), gl(m|n))$ form Howe
dual reductive pairs on $S(\C^d\otimes\C^{m|n})$ and on
${}^{Sp}H$, respectively. Therefore we have a decomposition of
modules
\begin{align*}
&\C[\x,\etabf]\cong\sum_{\la}V^\la_{Sp(d)}\otimes V^{\la'}_{osp(2m|2n)},\\
&{}^{Sp}H\cong\sum_{\la}V^\la_{Sp(d)}\otimes V^{\la''}_{m|n},
\end{align*}
where $\la$ is summed over a set of irreducible $Sp(d)$-highest
weights.  Here $\la'$ and $\la''$ are certain non-isomorphic
irreducible $osp(2m|2n)$- and $gl(m|n)$-highest weights,
respectively. Furthermore the map
${}^{Sp}I\otimes{}^{Sp}H\rightarrow\C[\x,\etabf]$ given by
multiplication is surjective and we have, for each $\la$,
$V^{\la'}_{osp(2m|2n)}={}^{Sp}IV^{\la''}_{m|n}$.
\end{thm}

The proofs of \thmref{O-duality} and \thmref{Sp-duality} are based
on the fact that the invariants of the classical group of the
corresponding dual pair in the endomorphism ring of
$S(\C^{d}\otimes\C^{m|n})$ are generated by quadratic invariants.
Although in \cite{H1} it is shown that the pairs
$(O(d),spo(2m|2n))$ and $(Sp(d),osp(2m|2n))$ are indeed dual pairs
on $S(\C^d\otimes\C^{m|n})$, the explicit decomposition of
$S(\C^d\otimes\C^{m|n})$ was not given.  We will embark on this
task in \secref{jointhwv}.

We conclude this section by showing that the representations of
$spo(2m|2n)$ and $osp(2m|2n)$ that appear in \thmref{O-duality}
and \thmref{Sp-duality} are unitarizable. We first recall some
definitions.

Let $A$ be a superalgebra and ${}^\dagger$ an anti-linear map with
$(ab)^\dagger=b^\dagger a^\dagger$, for $a,b$ in $A$. We call
${}^\dagger$ an anti-linear {\em anti-involution} if
$(a^\dagger)^\dagger=a$. Now let $A$ be a superalgebra equipped
with an anti-linear anti-involution ${}^\dagger$ and let $V$ be an
$A$-module. A Hermitian form $(\cdot|\cdot)$ on $V$ is said to be
{\em contravariant} if $(av|w)=(v|a^\dagger w)$, for $a\in A$ and
$v,w\in V$. If furthermore $(\cdot|\cdot)$ is positive-definite,
then $V$ is said to be a {\em unitarizable} $A$-module. We remark
here that we have defined the anti-involution and the
contravariant form without ``super signs''. It follows then that
any unitarizable module is completely reducible.

\begin{prop}\label{unitary}
The representations $V^{\la'}_{spo(2m|2n)}$ and
$V^{\la'}_{osp(2m|2n)}$ that occur in the decompositions of
$S(\C^d\otimes\C^{m|n})$ are unitarizable.
\end{prop}
\begin{proof} We need to construct a contravariant positive-definite
Hermitian form on $\C[\x,\etabf]$.  We proceed as follows.

First note that the space $\C[\x,\etabf]$ is an irreducible
representation of the direct sum of a Heisenberg algebra and a
Clifford superalgebra with generators mapped to $x_{i}^j$,
$\eta^j_k$, $\frac{\pa}{\pa x_{i}^j}$, $\frac{\pa}{\pa \eta^j_k}$,
for $i=1,\cdots,m$, $k=1,\cdots,n$ and $j=1,\cdots,d$, and $1$.
Identifying this superalgebra with its image we have an
anti-linear anti-involution given by
\begin{align*}
(x_{i}^j)^\dagger=\frac{\pa}{\pa x_{i}^j},\quad (\frac{\pa}{\pa
x_{i}^j})^\dagger=x_{i}^j,\quad
(\eta_{k}^j)^\dagger=\frac{\pa}{\pa \eta_{k}^j},\quad
(\frac{\pa}{\pa \eta_{k}^j})^\dagger=\eta_{k}^j,\quad 1^\dagger=1.
\end{align*}
 This gives
rise to a unique contravariant Hermitian form $(\cdot|\cdot)$ on
$\C[\x,\etabf]$ with $(1|1)=1$. Furthermore for any non-zero
monomial $f\in\C[\x,\etabf]$ we have $(f|f)>0$, and hence
$(\cdot|\cdot)$ is positive-definite. Therefore $\C[\x,\etabf]$,
as a representation of the Clifford superalgebra, is unitarizable.

Now it is easy to see, using \eqnref{glmn} along with the formulas
for ${}^{Sp}I$, ${}^{Sp}\Delta$ and ${}^OI$, ${}^O\Delta$ in this
section, that $osp(2m|2n)$ and $spo(2m|2n)$ are invariant under
the anti-involution ${}^\dagger$.  This implies that the
representations of $osp(2m|2n)$ and $spo(2m|2n)$ on
$\C[\x,\etabf]$ are unitarizable.
\end{proof}

\section{Joint highest weight vectors}\label{jointhwv}

In this section we will describe the explicit decomposition of the
space $S(\C^d\otimes\C^{m|n})$ under the joint actions of the
relevant dual pairs.  We will do so by explicitly finding a joint
highest weight vector for each irreducible component.

\subsection{The case of $(O,spo)$-duality}\label{hwvospo}
Consider the $(O(d),spo(2m|2n))$-duality on the space
$S(\C^d\otimes\C^{m|n})$.  Using the notation from
\secref{duality} we make the identification of
$S(\C^d\otimes\C^{m|n})$ with the polynomial superalgebra
$\C[\x,\etabf]$ so that the Lie algebra $so(d)$ and $spo(2m|2n)$
are identified with differential operators.

By \thmref{O-duality} we only need to find the decomposition of
the space of harmonic polynomials ${}^OH$ into irreducible
$O(d)\times gl(m|n)$-modules. By \thmref{glpmn-duality}
$\C[\x,\etabf]$ as a $gl(d)\times gl(m|n)$-module decomposes into
$\sum_{\la}V^\la_d\otimes V^\la_{m|n}$, where the summation is
over all partitions $\la$ with $l(\la)\le d$ and $\la_{m+1}\le n$.

Consider first the case when $m\ge\frac{d}{2}$. Take a diagram
$\la$ with $\la'_1+\la'_2\le d$ and let $v_{\la}$ be the
corresponding joint $gl(d)\times gl(m|n)$-highest weight vector in
$\C[\x,\etabf]$ of the form in \thmref{aux2} or
\thmref{glpmn-duality}. Note that in this case it is automatic
that $\la_{m+1}\le n$, as long as $n\ge 1$.

Here and further we use ${\bhf}$ to denote the $(m+n)$-tuple
$(\hf,\cdots,\hf;-\hf,\cdots,-\hf)$. That is, the first $m$
entries are $\hf$, while the last $n$ entries are $-\hf$.

\begin{prop}\label{O-harmonics}
Suppose that $n\ge 1$ and $m\ge\frac{d}{2}$. The vector $v_\la$ is
${}^O\Delta$-harmonic of $O(d)$-weight corresponding to the
diagram $\la$.  Therefore
\begin{equation*}
{}^OH\cong\sum_{\la}V^\la_{O(d)}\otimes V^{\la+d{\bf\hf}}_{m|n},
\end{equation*}
where $\la$ ranges over all diagrams with $\la'_1+\la'_2\le d$.
Here the weight $\la+d\bhf$ denotes the sum of the
$gl(m|n)$-weight corresponding to the Young diagram $\la$ with the
$(m+n)$-tuple $d\bhf$.
\end{prop}

\begin{proof}
Note that by our choice of the Borel subalgebra of $so(d)$, it is
automatic that $v_\la$ is an $O(d)\times gl(m|n)$-highest weight
vector.  (In fact this is true for any $\la$.)  Thus in order to
show that $v_\la$ is ${}^O\Delta$-harmonic it is enough to show
that it is annihilated by ${}^O\Delta^{xx}_{11}=\sum_{j=1}^d
\frac{\partial}{\partial x^j_1} \frac{\partial}{\partial
x^{d+1-j}_1}$. This is because $v_{\la}$ is already annihilated by
the nilpotent radical of the Borel subalgebra of $gl(m|n)$, which
together with ${}^O\Delta^{xx}_{11}{}$ generates the nilpotent
radical of the Borel subalgebra of $spo(2m|2n)$.  Also note that
if $\la'_1\le m$, then the joint highest weight vector is the
usual joint highest weight vector in the classical $O(d)\times
sp(2m)$-duality and hence is killed by ${}^O\Delta^{xx}_{11}$
\cite{H2}. So we may assume that $\la'_1>m$. In this case in order
to show that ${}^O\Delta^{xx}_{11}v_\la=0$, we consider the
classical $O(d)\times sp(2m+2n)$-duality.  Here the joint highest
weight vector $w_\la$ of $O(d)\times sp(2m+2n)$ is a product of
determinants of the form
$\Delta_{\la'_1}\Delta_{\la'_2}\cdots\Delta_{\la'_{\la_1}}$, with
only $\la'_1$ exceeding $m$. That is, only in $\Delta_{\la'_1}$
can we possibly have variables of the form $x_{m+i}^j$ with
$i=1,\cdots,n$.  From the duality in the classical case we know
that $w_\la$ is harmonic \cite{H2} and hence in particular
${}^O\Delta^{xx}_{11}w_\la=0$. Consider the first order
differential operators
$\Gamma_i=\sum_{j=1}^d\eta^{j}_1\frac{\pa}{\pa x^j_{m+i}}$, for
$i=1,\cdots,n$.  We see that ${}^O\Delta^{xx}_{11}$ commutes with
all $\Gamma_i$ and hence
\begin{equation*}
0=\Gamma_1\cdots\Gamma_{\la'_1-m}{}^O\Delta^{xx}_{11}w_\la=
{}^O\Delta^{xx}_{11}\Gamma_1\cdots\Gamma_{\la'_1-m}w_\la.
\end{equation*}
But $\Gamma_1\cdots\Gamma_{\la'_1-m}
w_\la=(-1)^{\la'_1-m-1}(\la'_1-m)!v_\la$ and hence
${}^O\Delta^{xx}_{11}v_\la=0$.

Finally the addition of $d{\bhf}$ to the $gl(m|n)$-highest weight
$\la$ is of course due to \eqnref{glmn}.
\end{proof}

As in \cite{H2} one shows that $v_\la$ indeed has $O(d)$-weight
corresponding to the Young diagram $\la$.  But as $\la$ ranges
over all partitions with $\la'_1+\la'_2\le d$ we conclude that the
$v_\la$'s generate the complete set of all finite-dimensional
irreducible $O(d)$-modules. Due to $O(d)\times gl(m|n)$-duality in
${}^OH$ we see that
$${}^OH\cong\sum_{\la}V^\la_{O(d)}\otimes V^{\la+d{\bhf}}_{m|n},$$
where $\la$ ranges over all Young diagrams with $\la'_1+\la'_2\le
d$.

\begin{prop}\label{aux3}
Suppose that $n\ge 1$ and $m\ge \frac{d}{2}$. Then as an
$O(d)\times spo(2m|2n)$-module we have the following
decomposition.
\begin{equation*}
S(\C^d\otimes\C^{m|n})\cong\sum_{\la}V^\la_{O(d)}\otimes
V^{\la+d\bf \hf}_{spo(2m|2n)},
\end{equation*}
where $\la$ ranges over all Young diagrams with $\la'_1+\la'_2\le
d$.
\end{prop}

Now consider the case when $m<\frac{d}{2}$.  For this case we
introduce new even variables so that the total number of even
variables is at least $\frac{d}{2}$.  Since the case when $d$ is
odd is analogous we assume for simplicity that $d$ is even and we
add new variables $x^j_{m+1},\cdots,x^j_{\frac{d}{2}}$,
$j=1,\cdots,d$, to the polynomial superalgebra $\C[\x,\etabf]$ and
denote the resulting superalgebra by $\C[\x',\etabf]$. That is, we
are considering the embedding $S(\C^d\otimes\C^{m|n})\subseteq
S(\C^d\otimes\C^{\frac{d}{2}|n})$. Without further mentioning we
adopt the convention of adding a $'$ to operators, vectors
etc.~when we are regarding them as over
$S(\C^d\otimes\C^{\frac{d}{2}|n})$. So for example we denote the
corresponding Laplacian of $\C[\x',\etabf]$ by
${}^O{\Delta'}^{xx}_{ij}$, $1\le i,j\le \frac{d}{2}$, etc.~ and
call ${}^O\Delta'$-harmonic an element $f\in\C[\x',\etabf]$ that
is annihilated by all these Laplacians. We note that
${}^O\Delta^{xx}_{is}={}^O{\Delta'}^{xx}_{is}$, for $1\le i,s\le
m$ etc. Furthermore ${}^O{\Delta'}^{xx}_{ij}$, with either $i$ or
$j$ not in $\{1,\cdots,m\}$, is a sum of second order differential
operators, each of them involving differentiation with respect to
some of the new variables $x^j_{m+1},\cdots,x^j_{\frac{d}{2}}$
that we have introduced. It follows that if
$f\in\C[\x,\etabf]\subseteq\C[\x',\etabf]$, then $f$ is
${}^O\Delta'$-harmonic if and only if $f$ is
${}^O\Delta$-harmonic. Thus ${}^OH\subseteq{}^OH'$.

Now in $\C[\x',\etabf]$ we know that the subspace of
${}^O\Delta'$-harmonics is
$${}^OH'=\sum_{\la}V^\la_{O(d)}\otimes
V^{\la+d\bhf}_{\frac{d}{2}|n},$$ where $\la$ ranges over all Young
diagrams with $\la'_1+\la'_2\le d$ by \propref{O-harmonics}.  Let
$v'_\la$ be a joint highest weight vector in ${}^OH'$ for the
component $V^\la_{O(d)}\otimes V^{\la+d\bhf}_{\frac{d}{2}|n}$.
Then if the first column exceeds $\frac{d}{2}$ we have up to a
scalar multiple
\begin{equation*}
v'_\la=\Delta_{1{\la'_1}}\Delta_{\la'_2}\cdots\Delta_{\la'_{\la_1}}.
\end{equation*}
Otherwise we have up to a scalar multiple
\begin{equation*}
v'_\la=\Delta_{{\la'_1}}\Delta_{\la'_2}\cdots\Delta_{\la'_{\la_1}}.
\end{equation*}
Suppose $\la$ is such a diagram with $\la_{m+1}>n$. In this case
the $n$-th column of $\la$ exceeds $m$ and hence $\Delta_{\la'_n}$
contains at least one row with entries consisting entirely of
newly introduced variables. Now by \thmref{glhwv} and
\thmref{O-duality} all the $O(d)$-highest weight vectors of
highest weight $\la$ in ${}^OH'$ are, up to a scalar multiple, of
the form $\Delta^T$, where $T$ runs over all
$(\frac{d}{2}|n)$-semi-standard tableaux. But then it is not hard
to see that one of the rows in some $\Delta^T_i$ must consist
entirely of newly introduced variables so that $\Delta^T$ reduces
to zero when setting $x^j_{m+1}=\cdots=x^j_{\frac{d}{2}}=0$. Since
${}^OH\subseteq{}^OH'$, this implies that there are no
$O(d)$-highest weight vectors of highest weight $\la$ with
$\la_{m+1}>n$ in ${}^OH$.

On the other hand if $\la'_{1}+\la'_2\le d$ and $\la_{m+1}\le n$,
it is quite easy to see, using \thmref{glpmn-duality}, that
$v_\la$ (that is the $O(d)\times gl(m|n)$-joint highest weight
vector in $\C[\x,\etabf]$) is annihilated by
${}^O\Delta^{xx}_{11}=\sum_{j=1}^d \frac{\partial}{\partial x^j_1}
\frac{\partial}{\partial x^{d+1-j}_1}$, and hence $v_\la$ is
indeed ${}^O\Delta$-harmonic. Combining the results of this
section we have proved the following.

\begin{thm}\label{aux1} We
have the following decomposition of ${}^OH$ as an $O(d)\times
gl(m|n)$-module:
\begin{equation*}
{}^OH\cong\sum_{\la}V^\la_{O(d)}\otimes V^{\la+d\bhf}_{m|n},
\end{equation*}
where $\la$ ranges over all diagrams with $\la'_1+\la'_2\le d$ and
$\la_{m+1}\le n$. Thus as an $O(d)\times spo(2m|2n)$-module we
have:
\begin{equation*}
S(\C^d\otimes\C^{m|n})\cong\sum_{\la}V^\la_{O(d)}\otimes
V^{\la+d\bf \hf}_{spo(2m|2n)},
\end{equation*}
where $\la$ ranges over all Young diagrams with $\la'_1+\la'_2\le
d$ and $\la_{m+1}\le n$. Here the labels of the
$spo(2m|2n)$-highest weight of $V^{\la+d\bhf}_{spo(2m|2n)}$ with
respect to the Dynkin diagram of \secref{ospirrep} is obtained by
applying \eqnref{spolabels} to the $gl(m|n)$-weight $\la+d\bhf$.
\end{thm}

\begin{proof}
The preceding discussion already shows that the theorem hold when
$n\ge 1$ and $m< \frac{d}{2}$. Since in the case when
$m\ge\frac{d}{2}$ and $n\ge 1$ the condition $\la_{m+1}\le n$ is
vacuous, the theorem is true in this case due to \propref{aux3}.
But of course the case $n=0$ is the well-known classical case for
which the conclusion of the theorem hold as well.
\end{proof}

\begin{rem}
Partial results on the decomposition of the space
$S(\C^d\otimes\C^{m|n})$ with respect to the joint action of $
O(d)\times spo(2m|2n)$ were obtained earlier by Nishiyama in
\cite{N2} by constructing certain $O(d)\times spo(2m|2n)$-joint
highest weight vectors in $S(\C^d\otimes\C^{m|n})$.  However, the
full set of such joint highest weight vectors (and hence the
complete decomposition) was not obtained in there.
\end{rem}

\subsection{The case of $(Sp,osp)$-duality} Now consider the
action of the dual pair $(Sp(d),osp(2m|2n))$ on the space
$S(\C^d\otimes\C^{m|n})$. The procedure is similar to that of
\secref{hwvospo}.

In view of \thmref{Sp-duality} we again only need to find the
decomposition of the space ${}^{Sp}H$ with respect to the joint
action of $Sp(d)\times gl(m|n)$. According to
\thmref{glpmn-duality} $\C[\x,\etabf]$ as a $gl(d)\times
gl(m|n)$-module decomposes into $\sum_{\la}V^\la_d\otimes
V^\la_{m|n}$, where the summation is over all partitions $\la$
with $l(\la)\le d$ and $\la_{m+1}\le n$.

Let us first consider the case when $m\ge\frac{d}{2}$.  We take a
Young diagram $\la$ with $l(\la)\le\frac{d}{2}$ so that the
condition $\la_{m+1}\le n$ here is automatic. We recall from
\secref{parameterization} that the finite-dimensional irreducible
representations of $Sp(d)$ are parameterized by diagrams with
length not exceeding $\frac{d}{2}$.

\begin{prop}\label{aux6} Suppose that $m\ge\frac{d}{2}$ and
let $\la$ be a diagram with $l(\la)\le \frac{d}{2}$. Let $v_\la\in
V^\la_d\otimes V^\la_{m|n}$ be a $gl(d)\times gl(m|n)$-joint
highest weight vector in $\C[\x,\etabf]$. Then $v_\la$ is
${}^{Sp}\Delta$-harmonic of $Sp(d)$-weight corresponding to the
diagram $\la$.  Therefore
\begin{equation*}
{}^{Sp}H\cong\sum_{\la}V^\la_{Sp(d)}\otimes
V^{\la+d{\bf\hf}}_{m|n},
\end{equation*}
where $\la$ ranges over all diagrams with $l(\la)\le \frac{d}{2}$.
Here $\la+d\bhf$ denotes the sum of the $gl(m|n)$-weight
corresponding to $\la$ and the $(m+n)$-tuple $d\bhf$.
\end{prop}

\begin{proof}
Since $m\ge\frac{d}{2}$ and $l(\la)\le\frac{d}{2}$, the
$gl(d)\times gl(m|n)$-joint highest weight vector is of the form
$v_\la=\prod_{i=1}^{\la_1}\Delta_{\la'_i}$, that is, only the $\x$
variables are involved. Since the Borel subalgebra of $sp(d)$ is
contained in the standard Borel subalgebra of $gl(d)$, $v_\la$ is
an $sp(d)\times gl(m|n)$-highest weight vector. We need to show
that it is ${}^{Sp}\Delta$-harmonic.  For this it is again
sufficient to show that $v_\la$ is annihilated by
${}^{Sp}\Delta^{xx}_{12}=\sum_{j=1}^{\frac{d}{2}}
\Big{(}\frac{\partial}{\partial x^j_1} \frac{\partial}{\partial
x^{d+1-j}_2}-\frac{\partial}{\partial x^{d+1-j}_1}
\frac{\partial}{\partial x^j_2}\Big{)}$. But this is clear by the
classical $Sp(d)\times so(2m)$-duality \cite{H2}, because the
formulas for the joint highest weight vector $v_\la$ and for the
Laplacian ${}^{Sp}\Delta^{xx}_{12}$ in the classical case are
identical with our formulas here. Now the proposition follows from
\thmref{Sp-duality} together with the fact that we have
constructed an $Sp(d)$-highest weight vector corresponding to
every finite-dimensional irreducible $Sp(d)$-module.
\end{proof}

We now consider the case $m<\frac{d}{2}$. In this case the
condition $\la_{m+1}\le n$ is not an empty condition. Here we can
apply the idea of \secref{O-duality} by inserting enough new even
variables $x_{m+1}^j,\cdots,x_{\frac{d}{2}}^j$ and consider the
$Sp(d)\times osp(d|2n)$-duality on the space
$S(\C^d\otimes\C^{\frac{d}{2}|n})$.  We identify
$S(\C^d\otimes\C^{\frac{d}{2}|n})$ with $\C[\x',\etabf]$ as before
and regard $\C[\x,\etabf]\subset\C[\x',\etabf]$. Again we will use
$'$ to distinguish elements in $\C[\x',\etabf]$ from elements in
$\C[\x,\etabf]$. As in \secref{O-duality} it is easy to see that
an element $f\in\C[x,\etabf]\subset\C[\x',\etabf]$ is
${}^{Sp}\Delta$-harmonic if and only if it is
${}^{Sp}\Delta'$-harmonic and therefore
${}^{Sp}H\subset{}^{Sp}H'$.

Now by \propref{aux6} we have
\begin{equation*}
{}^{Sp}H'=\sum_{\la}V_{Sp(d)}^\la\otimes
V^{\la+d\bhf}_{\frac{d}{2}|n},
\end{equation*}
where $\la$ is summed over all partitions of length $l(\la)\le
\frac{d}{2}$. A joint highest weight vector $v_\la$ is given by
$\prod_{i=1}^{\la_1}\Delta_{\la'_i}$ and hence by \thmref{glhwv},
the set of $\Delta^T$'s, where $T$ runs over all
$(\frac{d}{2}|n)$-semi-standard tableaux of shape $\la$, is a
basis for the space of $Sp(d)$-highest weight vectors of highest
weight $\la$ in ${}^{Sp}H'$.

Now suppose that $\la_{m+1}>n$.  Let $T$ be a
$(\frac{d}{2}|n)$-semi-standard tableau and
$\Delta^T=\prod_{i}\Delta^T_i$.  It is clear that in this case one
of the $\Delta^T_i$'s must contain a row consisting entirely of
newly introduced variables.  But then this means that, by setting
these newly introduced variables equal to zero, $\Delta^T$ is
zero.  This implies that in ${}^{Sp}H$ there are no
$Sp(d)$-highest weight vectors of highest weight $\la$, and hence
no $Sp(d)$-module of the form $V^\la_{Sp(d)}$ can occur in the
decomposition of ${}^{Sp}H$ with respect to the action of $Sp(d)$.

The above argument combined with \propref{aux6} gives the complete
description of the $Sp(d)\times osp(2m|2n)$-duality on the space
$S(\C^d\otimes\C^{m|n})$, which we summarize in the following
theorem.

\begin{thm}\label{Sposp-duality} We
have the following decomposition of ${}^{Sp}H$ as an $Sp(d)\times
gl(m|n)$-module:
\begin{equation*}
{}^{Sp}H\cong\sum_{\la}V^\la_{Sp(d)}\otimes V^{\la+d\bhf}_{m|n},
\end{equation*}
where $\la$ ranges over all diagrams with $l(\la)\le\frac{d}{2}$
and $\la_{m+1}\le n$. Thus as an $Sp(d)\times osp(2m|2n)$-module
we have:
\begin{equation*}
S(\C^d\otimes\C^{m|n})\cong\sum_{\la}V^\la_{Sp(d)}\otimes
V^{\la+d\bf \hf}_{osp(2m|2n)},
\end{equation*}
where $\la$ ranges over all Young diagrams with
$l(\la)\le\frac{d}{2}$ and $\la_{m+1}\le n$. Here the labels of
the $osp(2m|2n)$-highest weight of $V^{\la+d\bhf}_{osp(2m|2n)}$
with respect to the Dynkin diagram of \secref{ospirrep} is
obtained by applying \eqnref{osplabels} to the $gl(m|n)$-weight
$\la+d\bhf$.
\end{thm}

\begin{proof} Let $\la$ be a diagram with $l(\la)\le\frac{d}{2}$
and $\la_{m+1}\le n$. In view of the discussion above and
\propref{aux6} it remains to prove that in the case when
$m<\frac{d}{2}$, the $gl(d)\times gl(m|n)$-joint highest weight
vector in $\C[\x,\etabf]$ is indeed ${}^{Sp}\Delta$-harmonic. For
this it is enough to show that it is annihilated by
${}^{Sp}\Delta_{12}^{xx}=\sum_{j=1}^{\frac{d}{2}}\Big{(}
\frac{\partial}{\partial x^j_1} \frac{\partial}{\partial
x^{d+1-j}_2}-\frac{\partial}{\partial x^{d+1-j}_1}
\frac{\partial}{\partial x^j_2}\Big{)}$.  But this is easy to see
using the formula for such a joint highest weight vector given in
\thmref{glpmn-duality}.
\end{proof}

\section{Character formulas for irreducible unitarizable $spo(2m|2n)$- and
$osp(2m|2n)$-modules}\label{character}

In this section we give combinatorial character formulas for the
$spo(2m|2n)$- and $osp(2m|2n)$-representations that appear in the
decomposition of $S(\C^d\otimes\C^{m|n})$ of \secref{jointhwv}. We
shall need a result of Enright \cite{E} which we shall recall.
Before this we need some preparatory material.

Consider a Hermitian symmetric pair $(G,K)$, where $G$ is a real
classical simple Lie group. Let $\G$ and $\K$ denote the
corresponding complexified Lie algebras. Fix a Cartan subalgebra
$\h$ of $\K$ so that $\h$ is also a Cartan subalgebra of $\G$. Let
$\bb$ be a Borel subalgebra of $\G$ containing $\h$ so that ${\mf
q}=\K+\bb$ is a maximal parabolic subalgebra of $\G$ with abelian
radical ${\mf u}$. Hence as a vector space we have ${\mf
q}=\K\oplus{\mf u}$. Denote by $\Delta$ and $\Delta(\K)$ the root
systems of $(\G,\h)$ and $(\K,\h)$, respectively, and let
$\Delta_+$ be the set of positive roots determined by $\bb$.
Furthermore set $\Delta(\K)_+=\Delta_+\cap\Delta(\K)$ and let
$\rho$ and $\rho_\K$ denote the respective half sums of positive
roots. Also let $\Delta(\mf
u)=\{\alpha\in\Delta|\G_\alpha\subseteq{\mf u }\}$ and put
$\rho_{\mf u}=\hf\sum_{\alpha\in{\mf u}}\alpha$. Let $W$ and
$W(\K)$ denote the Weyl groups of $\G$ and $\K$, respectively.

Now to each $\la\in\h^*$ one can associate a subgroup $W_\la$ of
$W$.  Since we will need to explicitly compute $W_\la$ later on,
we will give a detailed description of it now.  The group $W_\la$
is the subgroup of $W$ generated by the reflections $s_\alpha$,
where $\alpha\in\Delta(\mf u)$ satisfying the following three
conditions \cite{E,DES}:
\begin{itemize}
\item[(i)] $<\la+\rho,\check{\alpha}>\in\N$.
\item[(ii)] If for some $\beta\in\Delta$ we have $(\la+\rho|\beta)=0$, then
$(\alpha|\beta)=0$.
\item[(iii)] If for some long root $\beta\in\Delta$ we have
$(\la+\rho|\beta)=0$, then $\alpha$ is a short root.
\end{itemize}
Associated to $W_\la$ one may define a root system $\Delta_\la$
consisting of the roots $\gamma\in\Delta$ such that $s_\gamma$
lies in $W_\la$. Now we set
$\Delta_{\la}(\K)=\Delta_\la\cap\Delta(\K)$, $\Delta_{\la
+}=\Delta_+\cap\Delta_\la$ and
$\Delta_\la(\K)_+=\Delta_\la(\K)\cap\Delta_{\la+}$. The group
$W_\la(\K)$ is defined to be the subgroup of $W_\la$ generated by
reflection along the roots lying in $\Delta_\la(\K)_+$. We have a
decomposition of the group $W_\la\cong W_\la(\K)\times W_\la^\K$,
where
\begin{equation}\label{Wklambda}
W_\la^\K=\{w\in W_\la|<w\rho,\check{\alpha}>\in\Z_+,\
\forall\alpha\in\Delta_{\la}(\K)_+\}.
\end{equation}

\begin{rem} Note
that our definition of $W_\la$ is actually the definition of
$W_{\la+\rho}$ in \cite{E} and \cite{DES}. \end{rem}

For $\mu\in\h^*$ being a $\Delta(\K)_+$-dominant integral weight
we denote the finite-dimensional irreducible $\K$-module of
highest weight $\mu$ by $V^\mu_{\K}$, as usual.

Now let $\la\in\h^*$ be a $\Delta(\K)_+$-dominant integral weight.
We may extend $V^\la_{\K}$ to a ${\mf q}$-module in the trivial
way and consider the induced representation $M^\la_{\G}$ of $\G$.
It is clear that $M^\la_{\G}$ contains a unique maximal submodule
and hence has a unique irreducible quotient, which is isomorphic
to the highest weight irreducible $\G$-module of highest weight
$\la$. We will denote this $\G$-module by $V^\la_{\G}$.

For $\xi\in\h^*$ with $<\xi,\check{\alpha}>\in\R$ for all
$\alpha\in\Delta(\K)$, we denote the unique
$\Delta(\K)_+$-dominant element in the $W(\K)$-orbit of $\xi$ by
$\overline{\xi}$.

We have the following character formula for an irreducible
unitarizable representation $V^\la_{\G}$.

\begin{thm}\cite{E,DES}\label{Enright} We have
\begin{equation*}
{\rm ch}V^\la_{\G}=\frac{e^{-\rho_{\mf u}}\sum_{w\in
W^\K_\la}(-1)^{l(w)}{\rm
ch}V_{\K}^{\overline{w(\la+\rho)}-\rho_\K}}
{\prod_{\alpha\in\Delta({\mf u })}(1-e^{-\alpha})},
\end{equation*}
where $l(w)$ is the length of $w$ in $W_\la$.
\end{thm}

\subsection{Character formula for $spo(2m|2n)$-modules}
It follows from \thmref{O-duality} in the case when $n=0$ that we
have the following identities of characters, for $d$ even and odd,
respectively.
\begin{align}
&(y_1\cdots y_m)^{\frac{d}{2}}
\prod_{i=1}^{\frac{d}{2}}\prod_{j=1}^m\frac{1}{(1-x_iy_j)(1-x^{-1}_iy_j)}
=\sum_{\la}{\rm ch}V_{O(d)}^\la{\rm
ch}V^{\la+d{\bf\hf}}_{sp(2m)},\ d{\rm\ even,}\label{evenchar}\allowdisplaybreaks\\
&(y_1\cdots y_m)^{\frac{d}{2}}
\prod_{i=1}^{\frac{d-1}{2}}\prod_{j=1}^m\frac{1}{(1-x_iy_j)
(1-x^{-1}_iy_j)(1-y_j)}=\sum_{\la}{\rm ch}V_{O(d)}^\la{\rm
ch}V^{\la+d{\bf\hf}}_{sp(2m)},\ d{\rm\ odd.}\label{oddchar}
\end{align}
Here $\la$ is summed over all partitions with $\la'_1+\la'_2\le d$
such that $l(\la)\le m$ and ${\bf\hf}$ stands for the $m$-tuple
$(\hf,\cdots,\hf)$. Let us write $\chi^{\la}_{O(d)}(\x)$ for the
character of $V^\la_{O(d)}$ to stress its dependence on the
variables $x_1,\cdots,x_{[\frac{d}{2}]}$. We will now apply
\thmref{Enright} to the Hermitian symmetric pair $(Sp(2m),U(m))$,
so that $\G=sp(2m)$ and $\K=gl(m)$. We may now rewrite ${\rm
ch}V^{\la+d{\bf\hf}}_{sp(2m)}$ in terms of Schur functions as
follows. Since $\rho_{\mf u}+d\bhf$ is $W(\K)$-invariant, we have
by \thmref{Enright}
\begin{equation*}
{\rm ch}V^{\la+d{\bf\hf}}_{sp(2m)}=(y_1\cdots
y_m)^{\frac{d}{2}}\frac{\sum_{w\in
W^\K_{\la+d\bhf}}(-1)^{l(w)}s_{\overline{w(\la+\rho_d)}-\rho_d}(\y)}
{\prod_{1\le i\le j\le m}(1-y_i y_j)},
\end{equation*}
where here and further $\rho_d=\rho+d\bhf$. Here
$W^{\K}_{\la+d\bhf}$ is the subset of the Weyl group of $sp(2m)$
defined by \eqnref{Wklambda}.

\begin{rem} As we now need to deal with $W_{\la+d\bhf}^{\K}$,
$W_{\la+d\bhf}$ and $W_{\la+d\bhf}(\K)$ for different $m$ at the
same time, we introduce a superscript $m$ in order to distinguish
them. So for example $W_{\la+d\bhf}^{\K,m}$ is the subset
$W_{\la+d\bhf}^\K$ of the Weyl group of $sp(2m)$.
\end{rem}

Combining this with \eqnref{evenchar} and \eqnref{oddchar},
respectively, we have for even and odd $d$ respectively
\begin{align}
\prod_{i=1}^{\frac{d}{2}}\prod_{j=1}^m\frac{1}{(1-x_iy_j)(1-x^{-1}_iy_j)}
=\sum_{\la}\chi&^\la_{O(d)}(\x) \frac{\sum_{w\in
W^{\K,m}_{\la+d\bhf}}(-1)^{l(w)}s_{\overline{w(\la+\rho_d)}-\rho_d}(\y)}
{\prod_{1\le
i\le j\le m}(1-y_i y_j)},\label{evenchar1}\allowdisplaybreaks\\
\prod_{i=1}^{\frac{d-1}{2}}\prod_{j=1}^m\frac{1}{(1-x_iy_j)
(1-x^{-1}_iy_j)(1-y_j)}&=\nonumber\\
\sum_{\la}&\chi^\la_{O(d)}(\x) \frac{\sum_{w\in
W^{\K,m}_{\la+d\bhf}}(-1)^{l(w)}s_{\overline{w(\la+\rho_d)}-\rho_d}(\y)}
{\prod_{1\le i\le j\le m}(1-y_i y_j)}.\label{oddchar1}
\end{align}
Here $\chi^\la_{O(d)}(\x)=\chi^{\bar{\la}}_{O(d)}(\x)$ if and only
if $\bar{\la}$ is obtained from $\la$ by replacing the first
column of $\la$ by a column of length $d-\la'_1$.  That is, the
corresponding representation of the Lie algebra $so(d)$ are
isomorphic.  Here and further we denote by $\bar{\la}$ the Young
diagram obtained from $\la$ by replacing its first column by
$d-\la'_1$ boxes.

In order to distinguish such representations at the level of
characters in the case when $d$ is odd let us take $-I\in
O(d)\setminus SO(d)$ and let $\epsilon$ denote the eigenvalue of
$-I$ so that we have $\epsilon^2=1$. We may then rewrite
\eqnref{oddchar1} as
\begin{align}
\prod_{i=1}^{\frac{d-1}{2}}\prod_{j=1}^m&\frac{1}{(1-\epsilon
x_iy_j) (1-\epsilon x^{-1}_iy_j)(1-\epsilon y_j)}=\nonumber\\
&\sum_{\la}\chi^\la_{O(d)}(\epsilon,\x) \frac{\sum_{w\in
W^{\K,m}_{\la+d\bhf}}(-1)^{l(w)}s_{\overline{w(\la+\rho_d)}-\rho_d}(\y)}
{\prod_{1\le i\le j\le m}(1-y_i y_j)}\label{oddchar2},
\end{align}
where now $\chi^\la_{O(d)}(\epsilon,\x)$ is a polynomial in $\x$
and $\epsilon$ such that when setting $\epsilon=1$, we obtain
$\chi^\la_{O(d)}(\x)$. Now it is easy to see that if $\la$ is a
Young diagram and $\chi^\la_{SO(d)}(\x)$ is the corresponding
$SO(d)$-character of $\chi^\la_{O(d)}(\epsilon,\x)$, then
$\chi^\la_{O(d)}(\epsilon,\x)=\epsilon^{|\la|}\chi^\la_{SO(d)}(\x)$,
where $|\la|$ is the size of $\la$.  Hence we have
$\chi^{\bar{\la}}_{O(d)}(\epsilon,\x)=\epsilon
\chi^\la_{O(d)}(\epsilon,\x)$.

The identites \eqnref{evenchar1} and \eqnref{oddchar2} will be our
starting point for a character formula for unitary
$spo(2m|2n)$-modules. We need the following lemma.

\begin{lem}\label{lem1}
Suppose that $f^\la(\y)$ and $g^\la(\y)$ are power series in the
variables $\y$.
\begin{itemize}
\item[(i)] Suppose that $d$ is odd and
\begin{equation}\label{O-independence}
\sum_{\la}f^\la(\y)\chi^\la_{O(d)}(\epsilon,\x)=\sum_\la
g^\la(\y)\chi^\la_{O(d)}(\epsilon,\x),
\end{equation}
where the summation is over the full set of irreducible
finite-dimensional characters of $O(d)$. Then
$f^\la(\y)=g^\la(\y)$, for all $\la$.
\item[(ii)] Suppose that $d$ is even and
\begin{equation*}
\sum_{\la}f^\la(\y)\chi^\la_{O(d)}(\x)=\sum_\la
g^\la(\y)\chi^\la_{O(d)}(\x),
\end{equation*}
where the summation is over the full set of irreducible
finite-dimensional characters of $O(d)$.  Then
$f^\la(\y)+f^{\bar{\la}}(\y)=g^\la(\y)+g^{\bar{\la}}(\y)$.
\end{itemize}
\end{lem}

\begin{proof}
We shall only show (i), i.e.~for $d$ odd, as the case of $d$ even
is analogous (in fact easier). The argument is similar to the that
of \cite{CL}.

We multiply the identity \eqnref{O-independence} by the Weyl
denominator $D$ of the Lie group $SO(d)$ and using the Weyl
character formula for $\chi^\la_{SO(d)}(\x)=\frac{\sum_{w\in
W}(-1)^{l(w)}e^{w(\la+\rho)}}{D}$ we obtain
\begin{equation}\label{auxx}
\sum_{\la}f^\la(\y)\epsilon^{|\la|}\sum_{w\in
W}(-1)^{l(w)}e^{w(\la+\rho)}=\sum_\la
g^\la(\y)\epsilon^{|\la|}\sum_{w\in W}(-1)^{l(w)}e^{w(\la+\rho)}.
\end{equation}
Now as $\la$ ranges over all integral dominant weights, $\la+\rho$
ranges over all regular integral dominant weights of $SO(d)$.
Hence if $\la\not=\mu$ as $SO(d)$-dominant weights, then the set
of weights $\{w(\la+\rho),w(\mu+\rho)|w\in W\}$ are all distinct.
Clearly two weights $\la$ and $\mu$ are equal as $SO(d)$-dominant
weights if and only if $\mu=\bar{\la}$. Thus looking at the
coefficient of $e^{\la+\rho}$ in \eqnref{auxx} we obtain
\begin{equation*}
\epsilon^{|\la|}f_\la(\y)e^{\la+\rho}+
\epsilon^{|\bar{\la}|}f_{\bar{\la}}(\y)e^{\bar{\la}+\rho} =
\epsilon^{|\la|}g_\la(\y)e^{\la+\rho}+
\epsilon^{|\bar{\la}|}g_{\bar{\la}}(\y)e^{\bar{\la}+\rho}.
\end{equation*}
Since $\epsilon^{|\la|}\epsilon^{|\bar{\la}|}=\epsilon$, we
conclude that $f^\la(\y)e^{\la+\rho}=g(\y)e^{\la+\rho}$ and hence
$f^\la(\y)=g^\la(\y)$.
\end{proof}

From the identities \eqnref{oddchar2} and \eqnref{evenchar1} by
using \lemref{lem1} we obtain the following results for every
$m\in\N$:

In the case when $d$ is odd
\begin{align*}
\sum_{w\in
W^{\K,m+1}_{\la+d\bhf}}(-1)^{l(w)}s_{\overline{w(\la+\rho_d)}-\rho_d}&(y_1,\cdots,y_m,0)
\\ =&\sum_{w\in
W^{\K,m}_{\la+d\bhf}}(-1)^{l(w)}s_{\overline{w(\la+\rho_d)}-\rho_d}(y_1,\cdots,y_m).
\end{align*}
In the case when $d$ is even
\begin{align*}
\sum_{w\in
W^{\K,m+1}_{\la+d\bhf}}(-1)^{l(w)}s_{\overline{w(\la+\rho_d)}-\rho_d}&(y_1,\cdots,y_m,0)
\\ &+\sum_{w\in
W^{\K,m+1}_{\bar{\la}+d\bhf}}(-1)^{l(w)}s_{\overline{w(\bar{\la}+\rho_d)}-\rho_d}(y_1,\cdots,y_m,0)\\
=\sum_{w\in
W^{\K,m}_{\la+d\bhf}}(-1)^{l(w)}s_{\overline{w(\la+\rho_d)}-\rho_d}&(y_1,\cdots,y_m)
\\&+ \sum_{w\in
W^{\K,m}_{\bar{\la}+d\bhf}}(-1)^{l(w)}s_{\overline{w(\bar{\la}+\rho_d)}-\rho_d}(y_1,\cdots,y_m).
\end{align*}
This allows us to define, in the case when $d$ is odd, an element
$S_{sp}^\la(y_1,y_2,\cdots)$ in the inverse limit of symmetric
polynomials, that is uniquely determined by the property that
\begin{equation*}
S^\la_{sp}(y_1,y_2,\cdots,y_m,0,0,\cdots)=\sum_{w\in
W^{\K,m}_{\la+d\bhf}}(-1)^{l(w)}s_{\overline{w(\la+\rho_d)}-\rho_d}(y_1,\cdots,y_m).
\end{equation*}
Similarly we may define an element
$S^\la_{sp}(y_1,y_2,\cdots)+S^{\bar{\la}}_{sp}(y_1,y_2,\cdots)$ in
the case when $d$ is even.

\begin{rem}
The elements $S_{sp}^\la(y_1,y_2,\cdots)$ and
$S^\la_{sp}(y_1,y_2,\cdots)+S^{\bar{\la}}_{sp}(y_1,y_2,\cdots)$
are in general infinite sums of symmetric functions and hence are
strictly speaking not symmetric functions.  However, in these
infinite sums there are only finitely many summands for any fixed
degree.
\end{rem}

We now take the limit as $m\rightarrow\infty$ in
\eqnref{evenchar1} and \eqnref{oddchar2} and obtain the following
identities, respectively.
\begin{align}
&\prod_{i=1}^{\frac{d}{2}}\prod_{j=1}^{\infty}\frac{1}{(1-x_iy_j)
(1-x^{-1}_iy_j)} =\sum_{\la}\chi^\la_{O(d)}(\x)
\frac{S^\la_{sp}(\y)} {\prod_{1\le
i\le j}(1-y_i y_j)},\label{evenchar22}\allowdisplaybreaks\\
&\prod_{i=1}^{\frac{d-1}{2}}\prod_{j=1}^\infty\frac{1}{(1-\epsilon
x_iy_j) (1-\epsilon x^{-1}_iy_j)(1-\epsilon y_j)}=
\sum_{\la}\chi^\la_{O(d)}(\epsilon,\x) \frac{S^\la_{sp}(\y)}
{\prod_{1\le i\le j}(1-y_i y_j)}\label{oddchar23},
\end{align}
where $\la$ is summed over all $O(d)$-highest weights and
$\y=(y_1,y_2,\cdots)$.

The identities \eqnref{evenchar22} and \eqnref{oddchar23} follow
from the fact that setting $y_{m+1}=y_{m+2}=\cdots=0$, they reduce
to identities \eqnref{evenchar1} and \eqnref{oddchar2},
respectively. Thus the left-hand sides and the right-hand sides of
\eqnref{evenchar22} and \eqnref{oddchar23} give rise to the same
elements in the ring of the symmetric functions, respectively.

Recall that $\omega$, the involution of the ring of symmetric
functions which sends the complete symmetric functions to the
elementary symmetric functions, is defined by
$\omega(\prod_{j\in\N}\frac{1}{1-w_j})=\prod_{j\in\N}(1+w_j)$ (see
for example \cite{M}). We can now apply $\omega$ partially to the
variables $y_{m+1},y_{m+2},\cdots$. After that we set the
variables $y_{m+n+1}=y_{m+n+2}=\cdots=0$ and we obtain the
following identities ($z_l=y_{m+l}$, for $l=1,\cdots,n$).
\begin{align}
\prod_{i=1}^{\frac{d}{2}}\prod_{j=1}^m\prod_{l=1}^n&\frac{(1+x_iz_l)
(1+x^{-1}_iz_l)} {(1-x_iy_j)
(1-x^{-1}_iy_j)}\nonumber\\
&=\sum_{\la}\chi^\la_{O(d)}(\x) \frac{HS^\la_{sp}(\y,{\bf
z})\prod_{i,l}(1+y_i z_l)} {\prod_{1\le i \le j\le m}(1-y_i
y_j)\prod_{1\le l<k\le n}(1-z_lz_k)},\label{evenchar3}
\allowdisplaybreaks\\
\prod_{i=1}^{\frac{d-1}{2}}
\prod_{j=1}^m\prod_{l=1}^n&\frac{(1+\epsilon x_iz_l)(1+\epsilon
x_i^{-1}z_l)(1+\epsilon z_l)} {(1-\epsilon x_iy_j)
(1-\epsilon x^{-1}_iy_j)(1-\epsilon y_j)}=\nonumber\\
&\sum_{\la}\chi^\la_{O(d)}(\epsilon,\x) \frac{HS^\la_{sp}(\y,{\bf
z})\prod_{i,l}(1+y_i z_l)} {\prod_{1\le i \le j\le m}(1-y_i
y_j)\prod_{1\le l<k\le n}(1-z_lz_k)}\label{oddchar3}.
\end{align}
\begin{rem} We note that $\omega(\prod_{1\le l\le k
}\frac{1}{1-z_lz_k})=\prod_{1\le l<k}\frac{1}{1-z_lz_k}$. This
follows from the following identities:
\begin{align*}
\prod_{1\le l\le k
}\frac{1}{1-z_lz_k}=\sum_{\la}s_\la(z_1,z_2,\cdots),\\
\prod_{1\le
l<k}\frac{1}{1-z_lz_k}=\sum_{\mu}s_\mu(z_1,z_2,\cdots),
\end{align*}
where $\la$ is summed over all partitions with even row lengths,
and $\mu$ is summed over all partitions with even column
lengths.
\end{rem}

Let us now explain the term $HS^\la_{sp}(\y,{\bf z})$. Since
setting the variables $y_{m+n+1}=y_{m+n+2}=\cdots=0$ the
expression $S^\la_{sp}(\y)$ reduces to a finite sum whose summands
are Schur polynomials with coefficients $\pm 1$, it follows that
applying the involution $\omega$ to it, we obtain a sum whose
summands consists of hook Schur functions with coefficients $\pm
1$. In fact if $S^\la_{sp}(y_1,y_2,\cdots)=\sum_{\mu}
\epsilon_{\mu}s_\mu(y_1,y_2,\cdots)$, where $\epsilon_\mu=\pm 1$,
then (cf.~\cite{CL})
\begin{equation*}
\omega(S^\la_{sp}(y_1,y_2,\cdots))= \sum_{\mu}\epsilon_\mu
HS_\mu(y_1,\cdots,y_m;z_1,z_2,\cdots),
\end{equation*}
where $HS_\mu(y_1,\cdots,y_m;z_1,z_2,\cdots)$ is the hook Schur
function of \cite{BR} in the variables $y_1,\cdots,y_m$ and
$z_1,z_2,\cdots$ corresponding to the partition $\mu$. Next
setting the variables $z_{n+1}=z_{n+2}=\cdots=0$ we get the hook
Schur polynomial associated to $\mu$, which we denote by
$HS_\mu(y_1,\cdots,y_m;z_1,\cdots,z_n)$. One property of hook
Schur polynomials is that
$HS_{\mu}(y_1,\cdots,y_m;z_1,\cdots,z_n)$ is non-zero if and only
if $\mu$ lies in the $(m|n)$-hook, i.e.~$\mu_{m+1}\le n$. So if
$S^\la_{sp}(\y)=\sum_{\mu}\epsilon_\mu s_\mu(\y)$, then by
$HS^\la_{sp}(\y,{\bf z})$ we mean the expression
$$HS^\la_{sp}(\y,{\bf z})=\sum_{\mu}\epsilon_\mu
HS_\mu(y_1,\cdots,y_m;z_1,\cdots,z_n).$$ Therefore $\la$ in
\eqnref{evenchar3} and \eqnref{oddchar3} is summed over all
$O(d)$-highest weights $\la$ such that $\la_{m+1}\le n$.

From \thmref{aux1}, \lemref{lem1}, and the identities
\eqnref{evenchar3} and \eqnref{oddchar3} we obtain the following
theorem.

\begin{thm}\label{char1} Let $\la$ be a diagram of
\thmref{aux1} and let $V^{\la+d{\bf\hf}}_{spo(2m|2n)}$ be the
irreducible $spo(2m|2n)$-module corresponding to $V_{O(d)}^\la$
under the Howe duality. Here $\bhf$ is the $(m+n)$-tuple
$(\hf,\cdots,\hf;-\hf,\cdots,-\hf)$.
\begin{itemize}
\item[(i)] If $d$ is odd, then
\begin{align*}
{\rm ch}V^{\la+d{\bf\hf}}_{spo(2m|2n)}=(\y{\bf
z^{-1}})^{\frac{d}{2}}\frac{HS_{sp}^\la(\y,{\bf
z})\prod_{i,l}(1+y_i z_l)} {\prod_{1\le i \le j\le m}(1-y_i
y_j)\prod_{1\le l<k\le n}(1-z_lz_k)}.
\end{align*}
\item[(ii)] If $d$ is even, then
\begin{align*}
{\rm ch}V^{\la+d{\bf\hf}}_{spo(2m|2n)}+{\rm
ch}V^{\bar{\la}+d{\bf\hf}}_{spo(2m|2n)}=(\y{\bf
z^{-1}})^{\frac{d}{2}}\frac{\big{(}HS_{sp}^\la(\y,{\bf
z})+HS_{sp}^{\bar{\la}}(\y,{\bf z})\big{)}\prod_{i,l}(1+y_i z_l)}
{\prod_{1\le i \le j\le m}(1-y_i y_j)\prod_{1\le l<k\le
n}(1-z_lz_k)}.
\end{align*}
\end{itemize}
Here $\y{\bf z^{-1}}$ stands for the product $y_1\cdots y_m
z_1^{-1}\cdots z_n^{-1}$.
\end{thm}

\begin{rem}
The expression $HS^\la_{sp}(\y,{\bf z})$ in general involves an
infinite number of hook Schur functions, so the computation of
these characters is a highly non-trivial task.  In order to have a
method to compute them, it is necessary to have an explicit
description of $W^{\K,m}_{\la+d\bhf}$. We will do this in
\secref{group}. From this we will then show in
\secref{consequences} that the coefficients of the monomials in a
character of a fixed degree can be computed by computing a finite
number of hook Schur functions.
\end{rem}

\subsection{Character formula for $osp(2m|2n)$-modules}
As the arguments in this case are very similar to the one given in
the previous section, we will only sketch them here.

It follows from \thmref{Sposp-duality} in the case when $n=0$ that
we have the following identity of characters.
\begin{align}
&(x_1\cdots x_m)^{\frac{d}{2}}
\prod_{i=1}^{\frac{d}{2}}\prod_{j=1}^m\frac{1}{(1-x_jy_i)(1-x_jy_i^{-1})}
=\sum_{\la}{\rm ch}V_{Sp(d)}^\la{\rm
ch}V^{\la+d{\bf\hf}}_{so(2m)}.\label{spchar}
\end{align}
Here $\la$ is summed over all partitions with $l(\la)\le {\rm min
}(\frac{d}{2},m)$ and ${\bf\hf}$ stands for the $m$-tuple
$(\hf,\cdots,\hf)$. Let us write $\chi^{\la}_{Sp(d)}(\y)$ for the
character of $V^\la_{Sp(d)}$ to stress its dependence on the
variables $y_1,\cdots,y_{\frac{d}{2}}$. We now apply
\thmref{Enright} to the Hermitian symmetric pair
$(SO^*(2m),U(m))$, so that we have $\G=so(2m)$ and $\K=gl(m)$. By
\thmref{Enright} we can then write ${\rm
ch}V^{\la+d{\bhf}}_{so(2m)}$ in terms of Schur functions as
\begin{equation*}
{\rm ch}V^{\la+d{\bf\hf}}_{so(2m)}=(x_1\cdots
x_m)^{\frac{d}{2}}\frac{\sum_{w\in
W^{\K,m}_{\la+d\bhf}}(-1)^{l(w)}s_{\overline{w(\la+\rho_d)}-\rho_d}(\x)}
{\prod_{i<j}(1-x_i x_j)}.
\end{equation*}
Here $W^{\K,m}_{\la+d\bhf}$ is a subset of the Weyl group of
$so(2m)$. Thus we have the following identity.
\begin{align}
&
\prod_{i=1}^{\frac{d}{2}}\prod_{j=1}^m\frac{1}{(1-x_jy_i)(1-x_jy_i^{-1})}
=\sum_{\la}\chi_{Sp(d)}^\la(\y)\frac{\sum_{w\in
W^{\K,m}_{\la+d\bhf}}(-1)^{l(w)}s_{\overline{w(\la+\rho_d)}-\rho_d}(\x)}
{\prod_{i<j}(1-x_i x_j)}.\label{spchar1}
\end{align}

Analogous to the proof of \lemref{lem1} one proves the following
lemma.

\begin{lem}\label{lem2}
Suppose that $f^\la(\x)$ and $g^\la(\x)$ are power series in the
variables $\x$ and suppose that
\begin{equation}\label{Sp-independence}
\sum_{\la}f^\la(\x)\chi^\la_{Sp(d)}(\y)=\sum_\la
g^\la(\y)\chi^\la_{Sp(d)}(\y),
\end{equation}
where the summation is over the full set of irreducible
finite-dimensional characters of $Sp(d)$. Then
$f^\la(\x)=g^\la(\x)$, for all $\la$.
\end{lem}

From \lemref{lem2} and the identity \eqnref{spchar1} it follows
that
\begin{align*}
\sum_{w\in
W^{\K,m+1}_{\la+d{\bhf}}}(-1)^{l(w)}s_{\overline{w(\la+\rho)}-\rho}&(x_1,\cdots,x_m,0)
=\\
&\sum_{w\in
W^{\K,m}_{\la+d\bhf}}(-1)^{l(w)}s_{\overline{w(\la+\rho_d)}-\rho_d}(x_1,\cdots,x_m),
\end{align*}
which then allows us to define an element
$S_{so}^\la(x_1,x_2,\cdots)$ in the inverse limit of symmetric
polynomials, uniquely determined by the property that
\begin{equation*}
S^\la_{so}(x_1,x_2,\cdots,x_m,0,0,\cdots)=\sum_{w\in
W^{\K,m}_{\la+d\bhf}}(-1)^{l(w)}s_{\overline{w(\la+\rho_d)}-\rho_d}(x_1,\cdots,x_m).
\end{equation*}

Taking the limit as $m\rightarrow\infty$ \eqnref{spchar1} and
\lemref{lem2} imply the following identity.
\begin{align}
&
\prod_{i=1}^{\frac{d}{2}}\prod_{j=1}^\infty\frac{1}{(1-x_jy_i)(1-x_jy_i^{-1})}
=\sum_{\la}\chi_{Sp(d)}^\la(\y)\frac{S^\la_{so}(x_1,x_2,\cdots)}
{\prod_{i<j}(1-x_i x_j)}.\label{spchar2}
\end{align}
We apply to \eqnref{spchar2} the involution of symmetric functions
$\omega$ partially to the variables $x_{m+1},x_{m+2},\cdots$, then
set the variables $z_{n+1}=z_{n+2}=\cdots=0$. We arrive at the
following identity ($z_l=x_{m+l}$, for $l=1,2,\cdots$).
\begin{align}
\prod_{i=1}^{\frac{d}{2}}\prod_{j=1}^m\prod_{l=1}^n&\frac{(1+y_iz_l)
(1+y^{-1}_iz_l)} {(1-y_ix_j)
(1-y^{-1}_ix_j)}\nonumber\\
&=\sum_{\la}\chi^\la_{Sp(d)}(\y) \frac{HS_{so}^{\la}(\x;{\bf
z})\prod_{i,l}(1+x_i z_l)}{\prod_{1\le i < j\le m}(1-x_i
x_j)\prod_{1\le l\le k\le n}(1-z_lz_k)},
\end{align}
where $HS^\la_{so}(x_1,\cdots,x_m;z_1,\cdots,z_n)$ is obtained by
applying the involution $\omega$ to $S^\la_{so}$ and setting the
variables $z_{n+1}=z_{n+2}=\cdots=0$.  As before it is also a sum
whose summands consist of hook Schur polynomials with coefficients
$\pm 1 $. By \thmref{aux3} and \lemref{lem2} we then obtain the
following theorem.

\begin{thm}\label{ospchar}
Let $\la$ be a diagram of \thmref{Sposp-duality} and let
$V^{\la+d{\bf\hf}}_{osp(2m|2n)}$ be the irreducible
$osp(2m|2n)$-module corresponding to $V_{Sp(d)}^\la$ under the
Howe duality. Here $\bhf$ is the $(m+n)$-tuple
$(\hf,\cdots,\hf;-\hf,\cdots,-\hf)$. Then
\begin{align*}
{\rm ch}V^{\la+d{\bf\hf}}_{osp(2m|2n)}={({\bf
xz^{-1}})}^{\frac{d}{2}}\frac{HS^{\la}_{so}(\x;{\bf
z})\prod_{i,l}(1+y_i z_l)} {\prod_{1\le i<j\le m}(1-x_i
x_j)\prod_{1\le l\le k\le n}(1-z_lz_k)},
\end{align*}
where ${\bf xz^{-1}}$ denotes the product $x_1x_2\cdots
x_mz^{-1}_1z^{-1}_2\cdots z^{-1}_n$.
\end{thm}

\begin{rem}
We actually have Howe dualities of the dual pairs
$(O(d),\G(C_{\infty}))$ and $(Sp(d),\G(D_{\infty}))$ on the space
$S(\C^d\otimes\C^\infty)$. Here the infinite-dimensional Lie
algebras $\G(C_\infty)$ and $\G(D_\infty)$ are Kac-Moody algebras
corresponding to the infinite affine matrices $C_\infty$ and
$D_\infty$, respectively \cite{K2}. From these dualities one can
show that, using similar arguments as we have given here, the
corresponding characters of those irreducible representations of
$\G(C_\infty)$- and $\G(D_\infty)$-modules are given by certain
infinite sums of symmetric functions. Applying the involution
$\omega$ to these characters one obtains the characters for our
$spo(2m|sn)$- and $osp(2m|2n)$-modules. Thus the characters of the
representations of $\G(C_{\infty})$ (respectively $\G(D_\infty)$)
that appear in these dualities determine the characters of the
representations $spo(2m|2n)$ (respectively $osp(2m|2n)$).
\end{rem}

\section{The Group $W_{\la+d\bhf}$}\label{group}

Throughout this section $\la=(\la_1,\la_2,\cdots,\la_{s})$ is a
partition of non-negative integers of length $s\le d$. We shall
describe the groups $W^m_{\la+d\bhf}$ and $W^m_{\la+d\bhf}(\K)$
for the Hermitian symmetric pairs $(Sp(2m),U(m))$ and
$(SO^*(2m),U(m))$.

Recall that the group $W^m_{\la+d{\bhf}}$ is defined as the
subgroup of the Weyl group of $sp(2m)$ or $so(2m)$ generated by
reflections corresponding to $\alpha\in\Delta({\mf u})$ satisfying
conditions (i), (ii) and (iii) given in \secref{character}.  We
will simply refer to them as conditions (i), (ii) and (iii) in
what follows.

\subsection{The case of $O(d)\times sp(2m)$-duality for $d$ even}
In the case when $d$ is even $W^m_{\la+d{\bhf}}$ is the subgroup
of the Weyl group of $sp(2m)$, which is isomorphic to the sign
permutation group $S_m\ltimes\Z^m_2$. The positive roots
$\Delta_+$ of $sp(2m)$ are generated by the simple roots
$-2\epsilon_1,\epsilon_1-\epsilon_2,
\epsilon_2-\epsilon_3,\cdots,\epsilon_{m-1}-\epsilon_m$.  We have
$\rho=-\epsilon_1-2\epsilon_2-\cdots-m\epsilon_m$, which we write
as
\begin{equation*}
\rho=(-1,-2,\cdots,-m).
\end{equation*}
We have the condition that $\la'_1+\la'_2\le d$. Now
$\Delta_+(\K)$ is generated by the simple roots
$\epsilon_1-\epsilon_2,\cdots,\epsilon_{m-1}-\epsilon_m$, while
$\Delta({\mf u})$ consists of roots of the form
$-\epsilon_i-\epsilon_j$, $1\le i\le j\le m$.

Let us first consider the case $s=\frac{d}{2}$. In this case
\begin{equation*}
\la+d{\bhf}+\rho=(\la_1+\frac{d}{2}-1,\la_2+\frac{d}{2}-2,\cdots,
\la_{\frac{d}{2}-1}+1,\la_{\frac{d}{2}},
-1,-2,\cdots,-m+\frac{d}{2}).
\end{equation*}
We see that $\la+d{\bhf}+\rho$ has no zero coefficient, and hence
condition (iii) is vacuous. It follows that for each
$i=1,\cdots,\frac{d}{2}$ with $m\ge\la_i+d-i$ we have
\begin{equation}\label{auxx2}
(\la+{d\bhf}+\rho,-\epsilon_i-\epsilon_{\la_i+d-i})=0.
\end{equation}
On the other hand if $m<\la_i+d-i$, $i=1\cdots,\frac{d}{2}$, we
have for all $t=1,\cdots,m$
\begin{equation}\label{auxxx2}
(\la+{d\bhf}+\rho,-\epsilon_i-\epsilon_t)<0.
\end{equation}
This implies by condition (ii) that if
$\alpha=-\epsilon_k-\epsilon_l$ is such that $s_\alpha\in
W^m_{\la+d{\bhf}}$, then neither $k$ nor $l$ can be in the index
set $J=\{1,\cdots,\frac{d}{2},\la_1+d-1,
\la_2+d-2,\cdots,\la_{\frac{d}{2}}+\frac{d}{2}\}$. Let
$I^{0}=\{1,\cdots,m\}\setminus J$. Let
$\alpha=-\epsilon_k-\epsilon_l$ with $k,l\in I^{0}$. Clearly we
have
\begin{equation*}
<\la+{d\bhf}+\rho,\check{\alpha}>\in\N,
\end{equation*}
and hence condition (i) is satisfied for such an $\alpha$. This
implies that $W^m_{\la+{d\bhf}}$ is generated by the reflections
$s_\alpha$ with $\alpha=-\epsilon_k-\epsilon_l$, $k,l\in I^{0}$.
Hence $W^m_{\la+{d\bhf}}$ is the sign permutation group on the
index set $I^0$.  Therefore $W^m_{\la+{d\bhf}}(\K)$ is equal to
the permutation group of the index set $I^{0}$ and hence
$\Delta_{\la+d{\bhf}}(\K)_+$ consists of $\epsilon_k-\epsilon_l$
with $k<l$ and $k,l\in I^{0}$.

Next consider the case $s<\frac{d}{2}$. In this case we have
\begin{equation*}
\la+d{\bhf}+\rho=(\la_1+\frac{d}{2}-1,\cdots,
\la_{s}+\frac{d}{2}-s,\frac{d}{2}-s-1,\cdots,
\underbrace{0}_{\frac{d}{2}},-1,-2,\cdots,-m+\frac{d}{2}).
\end{equation*}
Since $(\la+d{\bhf}+\rho,2\epsilon_{\frac{d}{2}})=0$, condition
(iii) implies that if $\alpha\in\Delta({\mf u})$ is such that
$s_\alpha\in W^m_{\la+d{\bhf}}$, then $\alpha$ is a short root. As
in the previous case \eqnref{auxx2} and \eqnref{auxxx2} hold in
this case as well with $i=1,\cdots,s$. In addition we have for
$j=s+1,\cdots,\frac{d}{2}$
\begin{equation}\label{auxx5}
(\la+d{\bhf}+\rho,-\epsilon_j-\epsilon_{d-j})=0.
\end{equation}
Let $J=\{1,\cdots,d-s-1,\la_1+d-1,\cdots,\la_s+d-s\}$ and
$I^{-}=\{1,\cdots,m\}\setminus J$. Similarly as in the previous
case conditions (i) and (ii) now tell us that
$\alpha\in\Delta({\mf u})$ is such that $s_\alpha\in
W^m_{\la+d{\bhf}}$ if and only if $\alpha=-\epsilon_k-\epsilon_l$
with $k,l\in I^-$ and $k\not=l$. Clearly $W^m_{\la+d{\bhf}}$ is
equal to the even sign permutation group (i.e.~permutations with
an even number of sign changes) in the index set $I^-$. Therefore
$W^m_{\la+d{\bhf}}(\K)$ is the permutation group on the index set
$I^-$ and hence $\Delta_{\la+d{\bhf}}(\K)_+$ consists of
$\epsilon_k-\epsilon_l$ with $k<l$ and $k,l\in I^{-}$.

Finally consider the case when $s>\frac{d}{2}$. In this case we
have
\begin{align*}
\la+d{\bhf}+\rho=&(\la_1+\frac{d}{2}-1,\cdots,
\underbrace{\la_{d-s}-\frac{d}{2}+s}_{d-s},s-\frac{d}{2},\\
&\cdots,\underbrace{1}_{\frac{d}{2}},0,\cdots,
\underbrace{1+\frac{d}{2}-s}_{s},
\underbrace{-1+\frac{d}{2}-s}_{s+1},\cdots,-m+\frac{d}{2}).
\end{align*}
Then \eqnref{auxx2} and \eqnref{auxxx2} hold for $i=1,\cdots,d-s$
and we have in addition
\begin{align}
&(\la+d{\bhf}+\rho,-\epsilon_j-\epsilon_{d-j+2})=0,\quad
j=d-s+2,\cdots,\frac{d}{2},\label{auxx4}\\
&(\la+d{\bhf}+\rho,-2\epsilon_{\frac{d}{2}+1})=0.\nonumber
\end{align}
Let
$I^+=\{d-s+1,s+1,s+2,\cdots,m\}\setminus\{\la_1+d-1,\cdots,\la_{d-s}+d-(d-s)\}$.
Then $W^m_{\la+d{\bhf}}$ is generated by $s_\alpha$, where
$\alpha=-\epsilon_k-\epsilon_l$ with $k,l\in I^+$ and $k\not=l$.
This implies that $W^m_{\la+d{\bhf}}$ is the even sign permutation
group on the index set $I^+$ and hence $W^m_{\la+d{\bhf}}(\K)$ is
the permutation group on the index set $I^+$ and hence
$\Delta_{\la+d{\bhf}}(\K)_+$ consists of $\epsilon_k-\epsilon_l$
with $k<l$ and $k,l\in I^+$.

\subsection{The case of $O(d)\times sp(2m)$-duality for $d$ odd}

Suppose that $s=\frac{d+1}{2}$. We have
\begin{equation*}
\la+d{\bhf}+\rho=(\la_1+\frac{d}{2}-1,\cdots,
\la_{\frac{d-3}{2}}-\frac{3}{2},\frac{3}{2},
\underbrace{\hf}_{\frac{d+1}{2}},-\frac{3}{2},
\cdots,-m+\frac{d}{2}).
\end{equation*}
Therefore \eqnref{auxx2} and \eqnref{auxxx2} hold for
$i=1,\cdots,\frac{d-3}{2}$ and also
$$(\la+d{\bhf}+\rho,-\epsilon_{\frac{d-1}{2}}
-\epsilon_{\frac{d+3}{2}})=0.$$ Note that the coefficients of
$\la+d{\bhf}+\rho$ are all half integers.  Hence if
$\beta\in\Delta({\mf u})$ is a long root, then
\begin{equation*}
<\la+d{\bhf}+\rho,\check{\beta}>\in\hf+\Z.
\end{equation*}
Thus the long roots are eliminated from the consideration of
$W^m_{\la+d{\bhf}}$ by condition (i). Now let
$J=\{1,\cdots,\frac{d-1}{2},\frac{d+3}{2},
\la_1+d-1,\cdots,\la_{d-s}+d-(d-s)\}$ and let
$I^0=\{1,\cdots,m\}\setminus J$. Then $W^m_{\la+d{\bhf}}$ is the
even sign permutation group on the index set $I^0$ so that
$W^m_{\la+d{\bhf}}(\K)$ is the permutation group on $I^0$.

Now suppose that $s<\frac{d+1}{2}$ so that
\begin{align*}
\la+d{\bhf}+\rho&=\\(\la_1+&\frac{d}{2}-1,\cdots,
\la_{s}+\frac{d}{2}-s,\frac{d}{2}-s-1,\cdots, \hf,
\underbrace{-\hf}_{\frac{d+1}{2}},-\frac{3}{2},\cdots,-m+\frac{d}{2}).
\end{align*}
Thus \eqnref{auxx2} and \eqnref{auxxx2} hold for $i=1,\cdots,s$.
Furthermore \eqnref{auxx5} holds for $
j=s+1,\cdots,\frac{d-1}{2}$. As before the long roots in
$\Delta({\mf u})$ are eliminated from considerations of
$W^m_{\la+d{\bhf}}$.  Set
$J=\{1,\cdots,d-s-1,\la_1+d-1,\cdots,\la_s+d-s\}$ and let
$I^-=\{1,\cdots,m\}\setminus J$. Then $W^m_{\la+d{\bhf}}$ is the
even sign permutation group on the index set $I^-$ so that
$W^m_{\la+d{\bhf}}(\K)$ is the permutation group on $I^-$.

Finally consider the case when $s>\frac{d+1}{2}$ so that we have
\begin{align*}
\la+d{\bhf}+\rho=&(\la_1+\frac{d}{2}-1,\cdots,
\underbrace{\la_{d-s}-\frac{d}{2}+s}_{d-s},s-\frac{d}{2},\\
&\cdots,\underbrace{\hf}_{\frac{d+1}{2}},\underbrace{-\hf}_{\frac{d+3}{2}},\cdots,
\underbrace{1+\frac{d}{2}-s}_{s},
\underbrace{-1+\frac{d}{2}-s}_{s+1},\cdots,-m+\frac{d}{2}).
\end{align*}
Thus \eqnref{auxx2} and \eqnref{auxxx2} still hold with
$i=1,\cdots,(d-s)$ and also \eqnref{auxx4} holds with
$j=1,\cdots,\frac{d+1}{2}$. Again the long roots are eliminated
from the consideration of $W^m_{\la+d{\bhf}}$. Let
$J=\{1,\cdots,(d-s),(d-s)+2,\cdots,{\frac{d+1}{2}},
\la_1+d-1,\cdots,\la_{d-s}+d-(d-s)\}$ and let
$I^+=\{1,\cdots,m\}\setminus J$. Then $W^m_{\la+d{\bhf}}$ is the
even sign permutation group on the index set $I^+$ so that
$W^m_{\la+d{\bhf}}(\K)$ is the permutation group on $I^+$.

\subsection{The case of $Sp(d)\times so(2m)$-duality}

In this case $W^m_{\la+d{\bhf}}$ is a subgroup of the Weyl group
of $so(2m)$, which is isomorphic to the even sign permutation
group $S_m\ltimes\Z_2^{m-1}$. The positive roots $\Delta_+$ is
generated by the simple root
$-\epsilon_1-\epsilon_2,\epsilon_1-\epsilon_2,
\epsilon_2-\epsilon_3,\cdots,\epsilon_{m-1}-\epsilon_m$, and hence
$\rho=-\epsilon_2-2\epsilon_3\cdots-(m-1)\epsilon_{m}$, which we
write as
\begin{equation*}
\rho=(0,-1,-2,\cdots,-m+1).
\end{equation*}
Let $\la=(\la_1,\la_2,\cdots,\la_{\frac{d}{2}})$ be a partition so
that
\begin{equation*}
\la+{d\bhf}+\rho=(\la_1+{\frac{d}{2}},\la_2+{\frac{d}{2}}-1,\cdots,
\la_{\frac{d}{2}}+1,\underbrace{0}_{\frac{d}{2}+1},
-1,-2,\cdots,-m+\frac{d}{2}+1).
\end{equation*}
The set $\Delta({\mf u})$ consists of roots of the form
$-\epsilon_k-\epsilon_l$, with $k\not=l$.

We have in the case $m\ge\la_i+d-i+2$
\begin{equation}\label{auxx3}
(\la+d{\bhf}+\rho,-\epsilon_i-\epsilon_{\la_i+d-i+2})=0,\quad
i=1,\cdots,\frac{d}{2}.
\end{equation}
On the other hand if $m<\la_i+d-i+2$, then for every
$t=1,\cdots,m$ we have
\begin{equation}\label{auxxx3}
(\la+d{\bhf}+\rho,-\epsilon_i-\epsilon_t)<0.
\end{equation}
This implies by condition (ii) that if
$\alpha=-\epsilon_k-\epsilon_l$ with $s_{\alpha}\in
W^m_{\la+d{\bhf}}$, then $k,l$ cannot be one of the indices in
\eqnref{auxx3} and \eqnref{auxxx3}. On the other hand set
$J=\{1,\cdots,\frac{d}{2},\la_1+d+1,\la_2+d,\cdots,
\la_{\frac{d}{2}}+\frac{d}{2}+2\}$ and let
$I=\{1,2,\cdots,m\}\setminus J$. Clearly if
$\alpha=-\epsilon_k-\epsilon_l$ with $k,l\in I$, then
\begin{equation*}
(\la+d{\bhf}+\rho,-\epsilon_k-\epsilon_l)\in\N,
\end{equation*}
and so condition (i) is satisfied.  Of course here (iii) is
irrelevant, as $\Delta$ is simply-laced.  Thus $W^m_{\la+d{\bhf}}$
is equal to the even sign permutation group on the index set $I$
and hence $W^m_{\la+d{\bhf}}(\K)$ is the permutation group on $I$.

From our explicit description of $W^m_{\la+d{\bhf}}$ we have the
following.

\begin{prop}
Let $\la=(\la_1,\cdots,\la_s)$ and let $W^m_{\la+d{\bhf}}$ be the
corresponding group of either sign or even sign permutation group
on the index set $I\subseteq\{1,\cdots,m\}$. Write
$\mu=\la+d{\bhf}+\rho=(\mu_1,\mu_2,\cdots,\mu_m)$ for the
generalized partition with either all integral or half-integral
row lengths.
\begin{itemize}
\item[(i)] We have
$|\mu_i|\not=|\mu_j|$, for $i,j\in I$ with $i\not=j$.
\item[(ii)] For $w\in W^m_{\la+d{\bhf}}$ the rows of the
generalized composition $w(\la+d{\bhf}+\rho)$ are all of different
length.
\item[(iii)] In the case $s\le\frac{d}{2}$ we have $\mu_i\le 0$, for
$i\in I$.
\item[(iv)] In the case $s>\frac{d}{2}$ we have $\mu_i>0$ and
$i\in I$ if and only if $i=(d-s+1)$ and
$\mu_{d-s+1}=s-\frac{d}{2}$.
\item[(v)] For all $m\in\N$ we have $W_{\la+d{\bhf}}^{m}\subseteq
W_{\la+d{\bhf}}^{m+1}$.
\end{itemize}
\end{prop}

\section{Consequences for the character
formula}\label{consequences}

In this section we will use the result of \secref{group} to study
the character formulas of \secref{character}.  In \secref{group}
we gave a description of $W^k_{\la+d{\bhf}}$.  However, in the
character formula we actually need to have a description of
$W^{\K,k}_{\la+d{\bhf}}$.

Let $w\in W^k_{\la+d{\bhf}}$.  From \secref{group} we know that
$w$ is either a sign permutation or an even sign permutation on an
index set $I$.

Recall the decomposition
$W^k_{\la+d{\bhf}}=W^k_{\la+d{\bhf}}(\K)\times
W_{\la+d{\bhf}}^{\K,k}$. As $W^k_{\la+d{\bhf}}$ is the sign or the
even sign permutation group on $I$ and $W^k_{\la+d{\bhf}}(\K)$ is
the permutation group on $I$, it follows that the elements of
$W_{\la+d{\bhf}}^{\K,k}$ are in one-to-one correspondence with
either the sign or the even sign changes of the index set $I$.
This correspondence can be made explicit as follows. Let
$I=\{j_1<j_2<\cdots<j_t\}$ and set
$\rho=(\rho_1,\rho_2,\cdots,\rho_k)$ and
$\rho_I=(\rho_{j_1},\rho_{j_2},\cdots,\rho_{j_t})$. Let $\tau$ be
either a sign change or an even sign change of $I$. Let $\sigma$
be the unique permutation on $I$ which permutes the rows of the
generalized composition $\tau(\rho_I)$ so that
$\sigma\tau(\rho_I)$ is a generalized partition. Set
$w_{\tau}=\sigma\tau$, then $w_\tau\in W_{\la+d{\bhf}}^{\K,k}$ is
the element corresponding to $\tau$ under the above-mentioned
one-to-one correspondence. More explicitly, if $\tau$ changes the
signs of $\rho_I$ at the rows $i_1<i_2<\cdots<i_l$, then $\sigma$
is the permutation that moves $i_l$ to $j_1$, $i_{l-1}$ to $j_2$,
$\cdots$, $i_1$ to $j_l$.  After that the remaining indices
$j_{l+1},\cdots,j_t$ are then assigned from the indices
$I\setminus\{i_1,\cdots,i_l\}$ in increasing order.

We are now in a position to describe
$\overline{w(\la+d{\bhf}+\rho)}$ with $w\in
W_{\la+d{\bhf}}^{\K,k}$. For this it is convenient to identify
$W_{\la+d{\bhf}}^{\K,k}$ with either a sign change or an even sign
change of the index set $I=\{j_1<j_2<\cdots<j_t\}$. Set
\begin{equation*}
\la+d{\bhf}+\rho=\mu=(\mu_1,\mu_2,\cdots,\mu_k).
\end{equation*}
Let $w\in W^{\K,k}_{\la+d{\bhf}+\rho}$ and let $\tau_w$ be the
corresponding sign change. Let us suppose that $\tau_w$
corresponds to sign changes of the subset $I_w\subseteq I$.
Suppose that $I_w={i_1,i_2,\cdots,i_l}$. Then
\begin{equation*}
\tau_w(\la+d{\bhf}+\rho)=(\mu_1,\mu_2,\cdots,-\mu_{i_1},
\cdots,-\mu_{i_2},\cdots,-\mu_{i_l},\cdots).
\end{equation*}
That is, $\tau_w(\mu)$ is obtained from $\mu$ by replacing all the
rows indexed by $I_w$ with its negative. Set $\sigma_w$ equal to
the unique permutation on $\{1,2,\cdots,k\}$ that permutes the
rows of generalized composition $\tau_w(\mu)$ so that the
resulting is a generalized partition. We denote by
$\Lambda_w(\la+d{\bhf}+\rho)$ the partition
$\sigma_w\tau_w(\la+d{\bhf}+\rho)-\rho-d{\bhf}$.

The following proposition is an easy consequence of the
correspondence between sign changes of the index set $I$ and
$W^{\K,k}_{\la+d{\bhf}}$.

\begin{prop}\label{highest1} With the notation
introduced above we have
\begin{equation*}
\Lambda_w(\la+d{\bhf}+\rho)= \overline{w(\la+d{\bhf}+\rho)}
-\rho-d{\bhf}.
\end{equation*}
\end{prop}

Using \propref{highest1} we can now prove the following corollary
for the characters of $V_{spo(2m|2n)}^{\la+d{\bhf}}$ and
$V_{osp(2m|2n)}^{\la+d{\bhf}}$. Recall the character formulas
given in \thmref{char1} and \thmref{ospchar}. In these formulas
the expression $HS_{\overline{w(\la+d{\bhf}+\rho)} -\rho-d{\bhf}}$
is the hook Schur function associated to the partition
${\overline{w(\la+d{\bhf}+\rho)} -\rho-d{\bhf}}$. Due to
\propref{highest1} we will from now on write
$HS_{\Lambda_w(\la+d{\bhf}+\rho)}$ for
$HS_{\overline{w(\la+d{\bhf}+\rho)} -\rho-d{\bhf}}$.

The next corollary shows that in general the character formulas
involve an infinite sum of hook Schur functions.

\begin{cor} Fix a diagram $\la$ corresponding to $V^{\la+d{\bhf}}_{spo(2m|2n)}$
or $V^{\la+d{\bhf}}_{osp(2m|2n)}$ of length $l(\la)=s$.
\begin{itemize}
\item[(i)] Suppose that $m>s+2$ and $n>1$. Then $HS^\la_{so}$ and
$HS^\la_{sp}$ are infinite sums of non-zero hook Schur
polynomials.
\item[(ii)] Suppose that $s=m$ and $\la_s>n$. Then $HS^\la_{so}$ and
$HS^\la_{sp}$ are finite sums of non-zero hook Schur polynomials.
\end{itemize}
\end{cor}

\begin{proof}
Take any $k>\la_1+d$ so that both $k$ and $k-1$ lie in the index
set $I$ associated to $W^{k}_{\la+d{\bhf}}$ of $sp(2k)$ or
$so(2k)$. Let $w\in W^{\K,k}_{\la+d{\bhf}}$ correspond to
$\tau_w$, the even sign permutation that permutes the indices
$k-1$ and $k$ and changes the signs of them. It is then easy to
see that $\La_w(\la+d{\bhf}+\rho)$ is a partition with $2$ from
the $(s+3)$-rd row on.  Since $n>1$, it follows that the partition
associated to $\La_w(\la+d{\bhf}+\rho)$ lies in the $(m|n)$-hook
and thus its corresponding hook Schur polynomial is non-zero. This
proves (i).

To prove (ii) let $k>\la_1+d$ and consider any $w\in
W^{\K,k}_{\la+d{\bhf}}$ corresponding to to a sign change
involving $k$.  Then $\La_w(\la+d{\bhf}+\rho)$ is a partition with
the first $m+1$ rows exceeding $n$.  But then the corresponding
hook Schur polynomial is zero.
\end{proof}

Let $\C[[\y,{\bf z}]]$ denote the ring of power series in the
variables $\y$ and ${\bf z}$.  We have a natural filtration of
ideals determined by the leading term.
$$\C[[\y,{\bf z}]]=\F_0\supset\F_1\supset\F_2\supset\cdots\supset\F_l\supset\cdots.$$

The formulas of \thmref{char1} and \thmref{ospchar} involve in
general an infinite number of hook Schur polynomials.  However,
for a fixed monomial that appears in the character formula we can
use a finite number of hook Schur polynomials to compute its
coefficient.  This follows from the following proposition.

\begin{prop} Let $\la$ be a partition with $l(\la)=s$ and $\la'_1+\la'_2\le d$. Let
$k>d$.
\begin{itemize}
\item[(i)] If $s>\frac{d}{2}$ we let $l=2k+|\la|-2s-1$. If
$s\le\frac{d}{2}$ we let $l=2k+|\la|-d$. Then we have
\begin{equation*}
HS_{sp}^\la(\y,{\bf z})\equiv\sum_{w\in
W^{\K,k-1}_{\la+d\bhf}}(-1)^{l(w)}HS_{\Lambda_w(\la+d{\bhf}+\rho)}
(\y;{\bf z}) \quad({\rm mod}\ \F_{l}).
\end{equation*}
\item[(ii)] Let $l=2k+|\la|-d-1$. Then
\begin{equation*}
HS_{so}^\la(\x,{\bf z})\equiv\sum_{w\in
W^{\K,k-1}_{\la+d\bhf}}(-1)^{l(w)}HS_{\Lambda_w(\la+d{\bhf}+\rho)}
(\x;{\bf z}) \quad({\rm mod}\ \F_{l}).
\end{equation*}
\end{itemize}
\end{prop}

\begin{proof} The theorem follows rather easily from \propref{highest1}.
We will only prove (ii), as (i) is quite similar. We may assume
without loss of generality that $k\in I$. Consider $w\in
W^{\K,k}_{\la+d{\bhf}}$ such that $w\not\in
W^{\K,k-1}_{\la+d{\bhf}}$. This means that $\tau_w$ changes the
sign of $k$. We consider the partition
$\Lambda_w(\la+d{\bhf}+\rho)$. It is not hard to see that the size
of this diagram is at least $2k-d+\sum_{i=1}^s\la_i-1$. But this
means that the hook Schur polynomial determined by
$\Lambda_w(\la+d{\bhf}+\rho)$ contains only monomials of degree
$l=2k-d+\sum_{i=1}^s\la_i-1$. Thus the hook Schur polynomials
associated to $\La_w(\la+d{\bhf}+\rho)$ with $w\in
W^{\K,k-1}_{\la+d{\bhf}}$ contain all monomials of $HS_{so}^\la$
of degree less than or equal to $l-1$.
\end{proof}

We now compute the functions $HS^\la_{so}$ and $HS^\la_{sp}$
explicitly in the case of $\la=(0,0,\cdots,0)$, the trivial
partition.

Let us first consider the case of $HS^\la_{sp}$ with $\la$ being
the trivial partition. We will write in this case simply $HS_{sp}$
for $HS_{sp}^\la$. In this case $W^k_{\la+d{\bhf}}$ is the group
of the even sign permutations in the indices $d,\cdots,k$.  Let
$w\in W^{\K,k}_{\la+d{\bhf}}$ and let $\tau_w$ be the
corresponding sign changes.  Let us suppose that $\tau_w$ changes
signs at the following $l$ rows: $i_1<i_2<\cdots<i_{l-1}<i_l$.
Here $i_1\ge d$ and $l$ is an even non-negative integer. Then it
is not hard to see that $\Lambda_w(d{\bhf}+\rho)$ is the following
partition
\begin{equation}\label{aux66}
(i_l-d+1,i_{l-1}-d+2,\cdots,i_1-d+l,\underbrace{l,\cdots,l}_{i_1-1},
\underbrace{l-1,\cdots,l-1}_{i_2-i_1-1},\cdots,\underbrace{1,\cdots,1}_{i_l-i_{l-1}-1}
,0,\cdots).
\end{equation}
That is, the first $l$ entries are
$i_l-d+1,i_{l-1}-d+2,\cdots,i_1-d+l$, followed by $i_1-1$ entries
of $l$ etc.  The length of $\Lambda_w(d{\bhf}+\rho)$ is $i_l$. For
a sequence of positive integers $I=\{i_1<i_2<\cdots<i_l\}$ with
$i_1\ge d$ denote by $\mu_I$ the partition \eqnref{aux66}.
Furthermore we let $|I|$ denote $l+\sum_{j=1}^l i_j$.

\begin{prop}\label{aux9} We have
\begin{equation*}
HS_{sp}(y_1,\cdots,y_m;z_1,\cdots,z_n)=
\sum_{I}(-1)^{|I|}HS_{\mu_I}(y_1,\cdots,y_m;z_1,\cdots,z_n),
\end{equation*}
where the summation is over all tuples $I=(i_1,i_2,\cdots,i_l)$
with $l$ even and $d\le i_1<i_2<\cdots<i_l$ satisfying one of the
following conditions.
\begin{itemize}
\item[(i)] In the case when $n\ge m$ we have $l\le n$ and at
most $m$ of the $i_j$'s exceed $d+n-m-1$.
\item[(ii)] In the case when $m>n$ we have $l<m$.
If in addition we have $l+i_t-t+1\le m+1\le l+i_{t+1}-t-1$, for
some $t=0,1,\cdots,l-1$, then $l-t\le n$. (Here by definition
$i_0$=0.)
\end{itemize}
\end{prop}

\begin{proof}
First we note that if $w\in W^{\K,k}_{\la+d{\bhf}}$ and $\tau_w$
its corresponding sign changes at the rows
$i_1<i_2<\cdots<i_{l-1}<i_l$, then $(-1)^{l(w)}=(-1)^{|I|}$.

(i) Obviously if $l>n$, then the corresponding partition $\mu_I$
cannot lie inside the $(m,n)$-hook.  Thus the corresponding hook
Schur polynomial is zero. Also clearly if $l\le m$, then $\mu_I$
lies in the $(m,n)$-hook. Now suppose that $m<l\le n$ and we have
\begin{equation*}
i_1<i_2<\cdots<i_{l-m}<i_{l-m+1}<\cdots<i_l.
\end{equation*}
Then $\mu_I$ lies in the $(m,n)$-hook if and only if
$i_{l-m}-d+(m+1)\le n$, which happens if and only if $i_{l-m}\le
d-m+n-1$.

(ii) Clearly, if $l\ge m$, then $\mu_I$ does not lie in the
$(m,n)$-hook. Now if $l<m$ and $i_l<m$, then it is easy to see
that $\mu_I$ lies in the $(m,n)$-hook.  On the other hand if $m\le
i_l$, we let $t=0,1,\cdots,l-1$ be such that
\begin{equation*}
l+i_t-t<m+1\le l+i_{t+1}-t-1.
\end{equation*}
It follows from \eqnref{aux66} that the $(m+1)$-st row of $\mu_I$
is $l-t$.  Thus $\mu_I$ lies in the $(m,n)$-hook if and only if
$l-t\le n$.
\end{proof}

Consider now the case of $HS^\la_{so}$, where $\la$ is the trivial
partition. We will again write in this case simply $HS_{so}$ for
$HS_{so}^\la$. Here $W^k_{\la+d{\bhf}}$ is the group of the even
sign permutations in the indices $d/2+1,d+2,d+3,\cdots,k$. Let
$w\in W^{\K,k}_{\la+d{\bhf}}$ and let $\tau_w$ be the
corresponding sign changes, which changes signs at the following
$l$ rows: $i_1<i_2<\cdots<i_{l-1}<i_l$.

First suppose that $i_1\not=d/2+1$. In this case
$\Lambda_w(d{\bhf}+\rho)$ is the partition
\begin{equation}\label{aux7}
(i_l-d-1,i_{l-1}-d,\cdots,i_1-d+l-2,\underbrace{l,\cdots,l}_{i_1-1},
\underbrace{l-1,\cdots,l-1}_{i_2-i_1-1},\cdots,\underbrace{1,\cdots,1}_{i_l-i_{l-1}-1}
,0,\cdots).
\end{equation}
Now if $i_1=\frac{d}{2}+1$, then $\Lambda_w(d{\bhf}+\rho)$ is
\begin{align}\label{aux8}
(i_l-d-1,i_{l-1}-d,\cdots,&i_2-d+l-3,\underbrace{l-1,\cdots,l-1}_{i_2-1},\\
&\underbrace{l-2,\cdots,l-2}_{i_3-i_2-1},\cdots,\underbrace{1,\cdots,1}_{i_l-i_{l-1}-1}
,0,\cdots).\nonumber
\end{align}
Note that \eqnref{aux8} is just \eqnref{aux7} corresponding to the
sequence $i_2<i_3<\cdots<i_l$, with $d+2\le i_2$.

For a sequence of positive integers $J=\{i_1<i_2<\cdots<i_l\}$
with $i_1\ge d+2$ we let $\nu_J$ be the partition \eqnref{aux7}.
Let $|J|=l+\sum_{j=1}^l i_j$.

\begin{prop}\label{aux10} We have
\begin{equation*}
HS_{so}(y_1,\cdots,y_m;z_1,\cdots,z_n)= \sum_{J}(-1)^{|J|}
HS_{\nu_J}(y_1,\cdots,y_m;z_1,\cdots,z_n),
\end{equation*}
where the summation is over all tuples $J=(i_1,i_2,\cdots,i_l)$
with $d+2\le i_1<i_2<\cdots<i_l$ satisfying the following
conditions.
\begin{itemize}
\item[(i)] In the case when $n\ge m$ we have $l\le n$ and at
most $m$ of the $i_j$'s exceed $d+n-m+1$.
\item[(ii)] In the case when $m>n$ we have $l<m$.
If in addition we have $l+i_t-t+1\le m+1\le l+i_{t+1}-t-1$, for
some $t=1,\cdots,l-1$, then $l-t\le n$. (Here again $i_0$=0.)
\end{itemize}
\end{prop}

\begin{proof}
As the proof is analogous to that of \propref{aux9}, we omit it.
\end{proof}

The module $V_{spo(2m|2n)}^\la$ (respectively
$V_{osp(2m|2n)}^\la$) with $\la$ being the trivial partition is
the $O(d)$-invariants (respectively $Sp(d)$-invariants) inside
$S(\C^d\otimes\C^{m|n})$.  Thus our computations of $HS_{sp}$ and
$HS_{so}$ give character formulas of these invariants.  On the
other hand we can describe the invariants, denoted by
$S(\C^d\otimes\C^{m|n})^{O(d)}$ and
$S(\C^d\otimes\C^{m|n})^{Sp(d)}$, in the following different way.
Since $gl(m|n)$ commutes with $O(d)$ and $Sp(d)$,
$S(\C^d\otimes\C^{m|n})^{O(d)}$ and
$S(\C^d\otimes\C^{m|n})^{Sp(d)}$ are modules over $gl(m|n)$. We
have the following analogue of classical invariant theory.

\begin{prop}\label{aux5} We have the following isomorphisms of
$gl(m|n)$-modules
\begin{itemize}
\item[(i)]
$S(\C^d\otimes\C^{m|n})^{O(d)}\cong\sum_{\la}V_{m|n}^{\la}$, where
the summation is over all partitions $\la$ with even row lengths,
$l(\la)\le d$ and $\la_{m+1}\le n$.
\item[(ii)]
$S(\C^d\otimes\C^{m|n})^{Sp(d)}\cong\sum_{\mu}V_{m|n}^{\mu}$,
where the summation is over all partitions $\mu$ with even column
lengths, $l(\mu)\le d$ and $\mu_{m+1}\le n$.
\end{itemize}
\end{prop}

\begin{proof}
The proof is in the same spirit as the one in the classical case
given in \cite{H2}. The $(gl(d),gl(m|n))$-duality gives
\eqnref{glgl-duality} and hence taking the $O(d)$-invariants on
both sides of \eqnref{glgl-duality} gives
\begin{equation*}
S(\C^d\otimes\C^{m|n})^{O(d)}=\sum_{\la} (V_d^\la)^{O(d)}\otimes
V_{m|n}^\la.
\end{equation*}
But it is known that $V_d^\la$ has only $O(d)$-invariants if and
only if $\la$ is an even partition, i.e.~all rows have even
length.  Furthermore in this case the dimension of
$O(d)$-invariants in $V_d^\la$ equals $1$.  This proves (i).

For (ii) we note that $V_d^\la$ has $Sp(d)$-invariants if and only
if $\la$ has even columns, in which case the dimension of the
invariants is again $1$.
\end{proof}

As the character of the $gl(m|n)$-module $V_{m|n}^\la$ is given by
the hook Schur function associated to $\la$ we obtain the
following corollary.

\begin{cor}
As $gl(m|n)$-modules we have
\begin{align*}
&{\rm ch}S(\C^d\otimes\C^{m|n})^{O(d)}=\sum_{\la}HS_{\la}
(y_1,\cdots,y_m; z_1,\cdots,z_n),\\
&{\rm
ch}S(\C^d\otimes\C^{m|n})^{Sp(d)}=\sum_{\mu}HS_{\mu}(x_1,\cdots,x_m;
z_1,\cdots,z_n),
\end{align*}
where the summations over $\la$ and $\mu$ are as in
\propref{aux5}.
\end{cor}

From these two descriptions of the $O(d)$-invariants inside
$S(\C^d\otimes\C^{m|n})$, in the case when $d$ is odd, we have the
combinatorial identity
\begin{align*}
\sum_{\la}HS_{\la}(y_1,\cdots,y_m&; z_1,\cdots,z_n)=\\
\sum_{I}(-1)^{|I|}HS_{\mu_I}(y_1,&\cdots,y_m;z_1,\cdots,z_n)
\Big{(}\frac{\prod_{1\le i\le m,1\le l\le n}(1+y_i z_l)}
{\prod_{1\le i \le j\le m,1\le l<k\le n}(1-y_i
y_j)(1-z_lz_k)}\Big{)},
\end{align*}
where $\la$ is summed over all partitions with even row lengths,
$l(\la)\le d$ and $\la_{m+1}\le n$, and $I$ is summed over all $I$
as in \propref{aux9} with $\mu_I$ as in \eqnref{aux66}.

Similarly from the descriptions of the $Sp(d)$-invariants we have
\begin{align*}
\sum_{\mu}HS_{\mu}(x_1,\cdots,x_m&; z_1,\cdots,z_n)=\\
\sum_{J}(-1)^{|J|} HS_{\nu_J}(x_1,&\cdots,x_m;z_1,\cdots,z_n)
\Big{(}\frac{\prod_{1\le i\le m,1\le l\le n}(1+x_i z_l)}
{\prod_{1\le i < j\le m,1\le l\le k\le n}(1-x_i
x_j)(1-z_lz_k)}\Big{)},
\end{align*}
where $\mu$ is summed over all partitions with even column
lengths, $l(\mu)\le d$ and $\mu_{m+1}\le n$, and $J$ is summed
over all $J$ as in \propref{aux10} with $\nu_I$ as in
\eqnref{aux7}.

\section{Tensor product decomposition}\label{tensor}

As another application of \thmref{aux1} and \thmref{Sposp-duality}
we derive in this section formulas for the decomposition of tensor
products of two representations of either $spo(2m|2n)$ or
$osp(2m|2n)$ that appear in the decomposition of
$S(\C^d\otimes\C^{m|n})$.

We first recall two Howe dualities involving the dual pairs
$(O(d),so(2k))$ and $(Sp(d),sp(2k))$ on the space
$\Lambda(\C^d\otimes\C^k)$.

\begin{thm}\label{aux90} \cite{H2}
The pairs $(O(d),so(2k))$ and $(Sp(d),sp(2k))$ form dual pairs on
the space $\Lambda(\C^d\otimes\C^k)$. Furthermore with respect to
their joint actions we have the following decompositions:
\begin{align}
&\Lambda(\C^d\otimes\C^k)\cong\sum_{\la} V^\la_{O(d)}\otimes
V^{\la'-d\bhf}_{so(2k)},\\
&\Lambda(\C^d\otimes\C^k)\cong\sum_{\mu} V^\mu_{Sp(d)}\otimes
V^{\mu'-d\bhf}_{sp(2k)},
\end{align}
where in the first sum $\la$ is summed over all diagrams with
$l(\la)\le d$, $\la'_1+\la'_2\le d$ and $\la_1\le k$, while in the
second sum $\mu$ is summed over all diagrams with $l(\mu)\le{d/2}$
and $\mu_1\le k$.
\end{thm}

\begin{rem}
We regard $so(2k)\cong osp(2k|0)$ and $sp(2k)\cong spo(2k|0)$ and
hence the labellings of their highest weights are as in
\secref{ospirrep}.
\end{rem}

Consider for positive integers $d$ and $r$ the decompositions
$S(\C^d\otimes\C^{m|n})\cong\sum_{\mu}V^\mu_{O(d)}\otimes
V^{\mu+d\bhf}_{spo(2m|2n)}$ and
$S(\C^r\otimes\C^{m|n})\cong\sum_{\gamma}V^\gamma_{O(r)}\otimes
V^{\gamma+r\bhf}_{spo(2m|2n)}$.  We have
\begin{align*}
S(\C^d\otimes\C^{m|n})\otimes
S(\C^r\otimes\C^{m|n})&\cong\sum_{\mu}V^\mu_{O(d)}\otimes
V^{\mu+d\bhf}_{spo(2m|2n)}\otimes\sum_{\gamma}V^\gamma_{O(r)}\otimes
V^{\gamma+r\bhf}_{spo(2m|2n)}\\
&\cong\sum_{\mu,\gamma}\big{(}V^\mu_{O(d)}\otimes
V^\gamma_{O(r)}\big{)}\otimes
\big{(}V^{\mu+d\bhf}_{spo(2m|2n)}\otimes
V^{\gamma+r\bhf}_{spo(2m|2n)}\big{)}.
\end{align*}
Now writing $V^{\mu+d\bhf}_{spo(2m|2n)}\otimes
V^{\gamma+r\bhf}_{spo(2m|2n)}\cong\sum_\la c^{\mu\gamma}_\la
V^{\mu+(d+r){\bhf}}_{spo(2m|2n)}$ we have therefore
\begin{equation}\label{aux91}
S(\C^d\otimes\C^{m|n})\otimes
S(\C^r\otimes\C^{m|n})\cong\sum_{\la,\mu,\gamma}c^{\mu\gamma}_\la
\big{(}V^\mu_{O(d)}\otimes V^\gamma_{O(r)}\big{)}\otimes
V^{\mu+(d+r){\bhf}}_{spo(2m|2n)}.
\end{equation}

On the other hand we have
\begin{align*}
S(\C^d\otimes\C^{m|n})\otimes S(\C^r\otimes\C^{m|n})&\cong
S(\C^{d+r}\otimes\C^{m|n})\\
&\cong\sum_\la V^\la_{O(d+r)}\otimes
V^{\la+(d+r){\bhf}}_{spo(2m|2n)}.
\end{align*}

If we let $V^\la_{O(d+r)}=\sum_{\mu,\gamma}b^\la_{\mu\gamma}
V^\mu_{O(d)}\otimes V^\gamma_{O(r)}$, that is, we regard
$V^\la_{O(d+r)}$ as an $O(d)\times O(r)$-module in the obvious
way, then we have
\begin{equation}\label{aux92}
S(\C^d\otimes\C^{m|n})\otimes S(\C^r\otimes\C^{m|n})\cong
\sum_{\la,\gamma,\mu} b^\la_{\mu\gamma}\big{(}V^\mu_{O(d)}\otimes
V^\gamma_{O(r)}\big{)}\otimes V^{\la+(d+r){\bhf}}_{spo(2m|2n)}.
\end{equation}

Combining \eqnref{aux91} and \eqnref{aux92} we see that
$c_\la^{\mu\gamma}=b^\la_{\mu\gamma}$.

This connection between the branching coefficients and the tensor
product coefficients, which may be regarded as a special case of
Kudla's seesaw pairs \cite{Ku}, is of course known \cite{H2}.

Now the same argument applied to the first dual pair of
\thmref{aux90} tells us that
$b_{\mu\gamma}^\la=a^{\mu\gamma}_\la$, where
\begin{equation*}
V_{so(2k)}^{\mu'-d\bhf}\otimes
V_{so(2k)}^{\gamma'-r\bhf}\cong\sum_\la a^{\mu\gamma}_\la
V_{so(2k)}^{\la'-(d+r){\bhf}}.
\end{equation*}
Taking account the fact that the $O(d)$-, $O(r)$- and
$O(2d)$-modules that appear in the various decompositions may not
be identical we have proved the following theorem.

\begin{thm} Let $\mu$ and $\gamma$ be diagrams lying in the $(m|n)$-hook and
satisfying the conditions $\mu_1'+\mu_2'\le d$ and
$\gamma_1'+\gamma_2'\le r$. Let $V^{\mu+d\bhf}_{spo(2m|2n)}\otimes
V^{\gamma+r\bhf}_{spo(2m|2n)}\cong\sum_{\la}c^{\mu\gamma}_\la
V^{\la+(d+r){\bhf}}_{spo(2m|2n)}$.  Let $k\ge {\rm
max}(\mu_1,\gamma_1)$ and $V_{so(2k)}^{\mu'-d\bhf}\otimes
V_{so(2k)}^{\gamma'-r\bhf}\cong\sum_\la a^{\mu\gamma}_\la
V_{so(2k)}^{\la'-(d+r){\bhf}}$.  Then for $\la$ lying in the
$(m|n)$-hook with $\la_1'+\la_2'\le d+r$ we have
$c_\la^{\mu\gamma}=a_\la^{\mu\gamma}$.  Otherwise
$c_\la^{\mu\gamma}=0$.
\end{thm}

We can derive the following theorem for $osp(2m|2n)$-modules in a
completely analogous fashion.

\begin{thm} For $d$ and $r$ even let $\mu$ and $\gamma$ be diagrams lying in the
$(m|n)$-hook with $l(\mu)\le d/2$ and $l(\gamma)\le r/2$. Let
$V^{\mu+d\bhf}_{osp(2m|2n)}\otimes
V^{\gamma+r\bhf}_{osp(2m|2n)}\cong\sum_{\la}c^{\mu\gamma}_\la
V^{\la+(d+r){\bhf}}_{osp(2m|2n)}$.  Let $k\ge {\rm
max}(\mu_1,\gamma_1)$ and $V_{sp(2k)}^{\mu'-d\bhf}\otimes
V_{sp(2k)}^{\gamma'-r\bhf}\cong\sum_\la a^{\mu\gamma}_\la
V_{sp(2k)}^{\la'-(d+r){\bhf}}$.  Then for $\la$ lying in the
$(m|n)$-hook with $l(\la)\le (d+r)/2$ we have
$c_\la^{\mu\gamma}=a_\la^{\mu\gamma}$.  Otherwise
$c_\la^{\mu\gamma}=0$.
\end{thm}

\begin{rem}
Of course the computation of the coefficients $a_\la^{\mu\gamma}$
are in general rather difficult.  There are combinatorial
algorithms that in principle can be used to compute them. See for
example \cite{Ki} and \cite{Li} and references therein. It turns
out that the coefficients can be computed once the usual
Littlewood-Richardson coefficients (for the general linear group)
are known.  The precise formulas are given in \cite{KT}.
\end{rem}

\begin{rem}
The tensor product decompositions of the $spo(2m|2n)$-modules and
the $osp(2m|2n)$-modules that appear in the decomposition of
$S(\C^d\otimes\C^{m|n})$ are stable in the following sense. The
coefficients $c_\la^{\mu\gamma}$ are independent of $m$ and $n$
for $n\ge 1$ and $m\ge d/2$.  This follows from a minor
modification of our argument above.
\end{rem}

\begin{rem}
The above method for computing the tensor product decomposition
using Howe duality appears to be quite general and could have
further applications. For example, using the $gl(d)\times
gl(m|n)$-Howe duality of \secref{glglduality} one can derive
rather easily the fact that the multiplication rule of the Hook
Schur functions is the same as that of ordinary Schur functions.
This was derived earlier in \cite{Re} using purely combinatorial
methods.
\end{rem}

\bigskip

\noindent{\bf Acknowledgements.} The second author gratefully
acknowledges partial financial support from the National Science
Council of the R.O.C. He also wishes to thank the Department of
Mathematics at National Taiwan University for hospitality.

\end{document}